# Functional Data-Driven Quantile Model Averaging with Application to Cryptocurrencies

Wenchao Xu[1], Xinyu Zhang[2], Jeng-Min Chiou[3,4] and Yuying Sun[2]

[1]*Shanghai University of International Business and Economics*

[2]*Academy of Mathematics and Systems Science, Chinese Academy of Sciences*

[3]*National Taiwan University and* [4]*Academia Sinica*

*Abstract:* Given the high volatility and susceptibility to extreme events in the cryptocurrency market, forecasting tail risk is of paramount importance. Value-at-Risk (VaR), a quantile-based risk measure, is widely used for assessing tail risk and is central to monitoring financial market stability. In data-rich environments, functional data from various domains are employed to forecast conditional quantiles. However, the infinite-dimensional nature of functional data introduces uncertainty. This paper addresses this uncertainty problem by proposing a novel data-driven conditional quantile model averaging (MA) approach. With a set of candidate models varying by the number of components, MA assigns weights to each model determined by a K-fold cross-validation criterion. We prove the asymptotic optimality of the selected weights in terms of minimizing the excess final prediction error when all candidate models are misspecified. Additionally, when the true regression relationship belongs to the set of candidate models, we provide consistency results for the averaged estimators. Numerical studies indicate that, in most cases, the proposed method outperforms other model selection and averaging methods, particularly for extreme quantiles in cryptocurrency markets.

*Keywords:* Quantile prediction, Model averaging, Asymptotic optimality, Functional data.

## 1. Introduction

The cryptocurrency market has experienced remarkable growth in recent years, attracting immense interest from a diverse array of stakeholders including investors, regulators, policymakers, companies, central banks, and governments worldwide (Cong et al. 2023, Guo et al. 2024, Sockin and Xiong 2023). Notably, these cryptocurrencies exhibit substantial fluctuations in value and are prone to occurrences of extreme price movements. For instance, Bitcoin, the pioneering cryptocurrency introduced in 2009 based on the scheme proposed by Nakamoto et al. (2008), attained an unprecedented global market capitalization of $2.54 trillion on May 12, 2021, but subsequently experienced a precipitous decline of nearly 50% to $1.3 trillion as of June 27, 2021. Consequently, forecasting tail risk, which refers to the potential for infrequent yet severe adverse events, is of critical importance across various domains, particularly in investment and portfolio management (Ando et al. 2022, Batur and Choobineh 2021, Capponi and Rubtsov 2022).

The central question we address in this paper is how to improve the accuracy of tail risk forecasts in cryptocurrency markets, given their high volatility and susceptibility to extreme events. We capture tail risk using Value-at-Risk (VaR), a conditional quantile that summarizes information regarding the distribution of potential losses. A major difficulty in forecasting tail risk stems from the substantial model uncertainty associated with the tail part of a loss distribution (Kou and Peng 2016). Uncertainty may arise from limited data, unknown models, and model misspecification (Hansen 2007, Hansen and Racine 2012). Furthermore, in the data-rich era, novel data structures like functional data, such as high-frequency asset returns (Saart and Xia 2022), contain more information than conventional point-valued data. As a result, functional data have emerged as popular covariates for forecasting conditional quantiles (Cardot et al. 2005), yet their use introduces further model uncertainty that must be accounted for.

Functional data are often inherently infinite-dimensional and thus functional principal component analysis (FPCA) for dimension reduction is widely used in functional regression (Cai and Hall 2006, Li et al. 2010, Kato 2012, Müller and Stadtmüller 2005, Yao et al. 2005b). One crucial issue is selecting the number of principal components to be retained for the covariate, which is essentially a model selection problem that aims to find the best model among a set of candidate models. Various model selection criteria include the fraction of variance explained (FVE, Chen and Müller 2012), leave-one-curve-out cross-validation (Kato 2012, Yao et al. 2005b), the Akaike information criterion (AIC, Kato 2012, Müller and Stadtmüller 2005, Yao et al. 2005b), and the Bayesian information criterion (BIC, Kato 2012). However, model selection assumes a unique "correct" model, though multiple competitive models may fit equally well. Choosing one optimal model ignores alternative model information, potentially putting all inferential weight on a flawed model. Moreover, individual models frequently exhibit sensitivity to data perturbations (Breiman 1996, 2001, Yuan and Yang 2005). Instead, model averaging (or forecast combination) incorporates all available information by constructing a weighted average over all candidate models. To the best of our knowledge, forecast combination of quantile regressions (QR) models with functional data remains unexplored in literature.

We propose a novel quantile model averaging (MA) approach for scalar responses with functional covariates, where the conditional quantile is modeled as a linear function of the functional covariate in Cardot et al. (2005). This flexible approach is motivated by the need to address the challenges posed by the unique characteristics of infinite-dimensional data and model uncertainties. To handle functional data, whether sparsely or densely observed, we employ principal components analysis through conditional expectation (Yao et al. 2005a) to estimate functional principal component scores. Our MA method considers a set of candidate models with varying numbers of principal components and assigns data-driven weights by minimizing $K$-fold cross-validation prediction errors. We provide two novel theoretical justifications for MA. First, we show that when there are no correctly specified models (i.e., true models) in the set of candidate models, the proposed method is asymptotically optimal in the sense that its excess final prediction error is as small as that of the

# 1 INTRODUCTION

infeasible best possible prediction. Second, in the situation with true candidate models, we show that the averaged estimators for the model parameters are consistent.

Our Monte Carlo simulations demonstrate that the proposed MA achieves lower risk than the competing methods when all candidate models are misspecified. We apply our method to predict the quantile of the minimum hourly log-return of ten cryptocurrencies and analyze the related tail risks. Empirical results indicate the superiority of the proposed approach compared with other model averaging and selection methods for VaR forecasts. For instance, for a given probability level $\tau$ and the cryptocurrency Bitcoin Cash, our method outperforms other model averaging methods, yielding nearly a twofold improvement in performance. By using calibration tests, we show that the proposed MA approach outperforms other competing methods. Furthermore, model averaging outperforms the related model selection techniques, which highlights the merits of reducing model uncertainty and increasing the robustness of model averaging.

It is worth discussing some key references and outlining our contributions in relation to the most relevant literature.

1. To the best of our knowledge, Zhang et al. (2018) and Zhang and Zou (2020) are the only existing studies on model averaging for prediction with functional data. However, their papers focus on the conditional mean prediction setting where both the response and the covariate are treated as random functions within the framework of functional regression models. In contrast, our work considers the conditional quantile prediction with functional covariates. This setting, where only the covariates are infinite-dimensional functional data while the responses are scalar values, is widely encountered across numerous applied domains.

2. We establish the asymptotic optimality based on the excess final prediction error which removes the dominant term arising from the innovations. This result improves upon the existing asymptotic optimality in the sense of minimizing the final prediction error in the literature on QR model averaging (Lu and Su 2015, Wang et al. 2023). This is an important contribution to the model averaging for QR. Moreover, our work establishes a convergence rate for the asymptotic optimality, which has not established for either conditional mean or conditional quantile processes in existing literature (Cheng and Hansen 2015, Feng et al. 2024, Hansen and Racine 2012, Nelson et al. 2021).

3. We are the first to establish the consistency of the weighted average estimator and derive its rate of convergence for QR under the condition that at least one of the candidate models is correctly specified. In contrast, most existing studies focus solely on the optimality property without establishing consistency or convergence rates.

The remainder of the paper is organized as follows. Section 2 lists literature review about quantile prediction with functional data and model averaging. Section 3 presents the model and FPCA-based estimation and conditional quantile prediction. Section 4 proposes a model averaging procedure for conditional quantile prediction. Section 5 establishes the asymptotic properties of the resulting model averaging estimator. Section 6 conducts simulation studies to illustrate its finite sample performance. Section 7 applies the proposed method to cryptocurrencies. Section 8 concludes the paper. All conditions and mathematical proofs are listed in the Supplementary Material.

## 2. Literature Review

Our literature review focuses on two different research aspects. First, we review the development of quantile prediction with functional data. Second, we analyze the literature on model averaging.

### 2.1. Quantile prediction with functional data

QR has attracted increasing attention since the seminal work introduced by Koenker and Bassett (1978) due to its ability to provide richer information about the conditional distribution of a response variable compared to traditional mean regression (Koenker 2005). With recent technological advancements in data-rich environments, observations are increasingly recorded as functional data, either continuously over time intervals or intermittently at discrete time points. Functional data analysis has become increasingly important over the past two decades (Ramsay and Silverman 2005, Ferraty and Vieu 2006, Horváth and Kokoszka 2012, Zhang 2013, Hsing and Eubank 2015, Cuevas 2014, Morris 2015, Wang et al. 2016). These functional data have been utilized to model the conditional quantiles of scalar responses, leading to the development of functional quantile regression (FQR), which extends the standard QR framework to account for functional covariates.

FQR has led to extensive research on estimation methods, asymptotic properties, and applications. For example, Cardot et al. (2005) considered a smoothing splines-based approach to represent the functional covariates for the functional linear quantile regression (FLQR) model and established a convergence rate result. Kato (2012) studied FPCA-based estimation for the FLQR model and established minimax optimal convergence rates for regression estimation and prediction. Li et al. (2022) and Sang et al. (2022) considered statistical inference for the FLQR model based on FPCA and reproducing kernel Hilbert space framework, respectively. Ferraty et al. (2005) and Chen and Müller (2012) estimated the conditional quantile function by inverting the corresponding conditional distribution function. Yao et al. (2017) and Ma et al. (2019) studied a high-dimensional partially FLQR model that contains vector and functional-valued covariates.

A key characteristic of functional data is its inherent infinite-dimensionality (Chen and Müller 2012, Yao et al. 2005a). To address this, dimensionality reduction techniques are employed to represent the functional covariates using a finite number of components (Kato 2012, Yao et al. 2005b). The choice of the number of components retained can significantly impact the performance of conditional quantile prediction. Moreover, functional covariates can be observed densely or sparsely, further compounding the challenges associated with prediction tasks involving such data. Model selection methods provide one approach to selecting the number of components, although they may not yield optimal results in terms of prediction error and carry the risk of selecting an inferior model. To address these difficulties in quantile prediction with functional covariates, this paper proposes a model averaging prediction by minimizing the $K$-fold cross-validation weight choice criterion.

## 2.2. Model averaging

Model averaging is an alternative approach to model selection that involves combining a set of candidate models by taking a weighted average of them. This can potentially reduce risk relative to model selection (Liu 2015, Liu and Okui 2013, Magnus et al. 2010, Peng and Yang 2022, Yuan and Yang 2005). Hansen and Racine (2012) proposed jackknife (leave-one-out cross-validation) model averaging for least squares regression that selects the weights by minimizing a cross-validation criterion. Jackknife model averaging was further carried out on models with dependent data (Zhang et al. 2013). Cheng and Hansen (2015) extended the jackknife model averaging method to leave-$h$-out cross-validation criteria for forecasting in combination with factor-augmented regression, and Gao et al. (2016) further extended it to leave-subject-out cross-validation under a longitudinal data setting. Recently, the model averaging method has been developed for functional data analysis. For example, Zhang et al. (2018) developed a cross-validation model averaging estimator for functional linear regression in which both the response and the covariate are random functions. Zhang and Zou (2020) proposed a jackknife model averaging method for a generalized functional linear model.

The weight selection criteria used in the model averaging mentioned above are designed for the conditional mean process, which is inappropriate for quantile prediction. To the best of our knowledge, there are two strategies designed for QR frequentist model averaging in the literature.

Specifically, Lu and Su (2015) constructed Mallows-type and jackknife model averaging for QR, which were subsequently extended by Wang et al. (2023) through a two-step process involving model screening for high-dimensional QR. However, these works employ traditional point-valued data as predictors and do not address quantile prediction issues with functional covariates. Therefore, it is necessary to construct a weight selection criteria specific to quantile prediction problems involving functional covariates, which poses unique challenges due to their infinite-dimensional nature.

## 3. Model and Estimation

### 3.1. Model Set-up

Let $\{(Y_i, X_i)\}_{i=1}^n$ be an independent and identically distributed (i.i.d.) sample, where $Y_i$ is a scalar random variable and $X_i = \{X_i(t): t \in \mathcal{T}\}$ is a random function defined on a compact interval $\mathcal{T}$ of $\mathbb{R}$. Let $Q_\tau(X_i)$ denote the $\tau$th conditional quantile of $Y_i$ given $X_i$, where $\tau \in (0, 1)$ is the quantile index. An FLQR model implies that $Q_\tau(X_i)$ can be written as a linear functional of $X_i$, that is, there exist a scalar constant $a(\tau) \in \mathbb{R}$ and a slope function $b(\cdot, \tau) \in L_2(\mathcal{T})$ such that

$$Q_\tau(X_i) = a(\tau) + \int_{\mathcal{T}} b(t, \tau) X_i^c(t)\, dt, \tag{3.1}$$

where $L_2(\mathcal{T})$ is the set of square-integrable functions defined on $\mathcal{T}$, $X_i^c(t) = X_i(t) - \mu(t)$, and $\mu(t) = E\{X_i(t)\}$ is the mean function of $X_i$. For notational simplicity, we suppress the dependence of $a(\tau)$ and $b(\cdot, \tau)$ on $\tau$. As a result, we have the following FLQR model

$$Y_i = Q_\tau(X_i) + \varepsilon_i = a + \int_{\mathcal{T}} b(t) X_i^c(t)\, dt + \varepsilon_i,$$

where $\varepsilon_i \equiv Y_i - Q_\tau(X_i)$ satisfies the quantile restriction $P(\varepsilon_i \leq 0 | X_i) = \tau$ almost surely (a.s.).

Denote the covariance function of $X_i(t)$ by $G(s, t) = \text{Cov}\{X_i(s), X_i(t)\}$. Since $G$ is a symmetric and nonnegative definite function, Mercer's theorem (Hsing and Eubank 2015, Theorem 4.6.5) ensures that there exist a non-increasing sequence of eigenvalues $\kappa_1 \geq \kappa_2 \geq \cdots > 0$ and a sequence of corresponding eigenfunctions $\{\phi_j\}_{j=1}^{\infty}$, such that $G(s,t) = \sum_{j=1}^{\infty} \kappa_j \phi_j(s) \phi_j(t)$. Since $\{\phi_j\}_{j=1}^{\infty}$ is an orthonormal basis of $L_2(\mathcal{T})$, we have the following expansions in $L_2(\mathcal{T})$:

$$X_i^c(t) = \sum_{j=1}^{\infty} \xi_{ij} \phi_j(t) \quad \text{and} \quad b(t) = \sum_{j=1}^{\infty} b_j \phi_j(t),$$

where $b_j = \int_{\mathcal{T}} b(t) \phi_j(t)\, dt$ and $\xi_{ij} = \int_{\mathcal{T}} X_i^c(t) \phi_j(t)\, dt$. The $\xi_{ij}$'s are called functional principal component scores and satisfy $E(\xi_{ij}) = 0$, $E(\xi_{ij}^2) = \kappa_j$, and $E(\xi_{ij}\xi_{ik}) = 0$ for all $j \neq k$. The expansion for $X_i^c(t)$ is called the Karhunen-Loève expansion (Hsing and Eubank 2015, Theorem 7.3.5). Then, the model (3.1) can be transformed into a linear QR model with an infinite number of covariates:

$$Q_\tau(X_i) = a + \sum_{j=1}^{\infty} b_j \xi_{ij}. \tag{3.2}$$

A nonlinear ill-posed inverse problem would be encountered when trying to fit the model (3.2) with a finite number of observations; see Kato (2012) for a detailed discussion. To address this issue, we follow Kato (2012) and truncate the model (3.2) to a feasible linear QR model. It is common to keep only the first $J_n < \infty$ eigenbases, which leads to the following truncated model

$$Q_{\tau,J_n}(X_i) = a + \sum_{j=1}^{J_n} b_j \xi_{ij}, \tag{3.3}$$

where we allow $J_n$ to increase to infinity as $n \to \infty$ to reduce approximation errors. The regression parameters $a, b_1, \ldots, b_{J_n}$ in model (3.3) need to be estimated.

### 3.2. FPCA-based Quantile Prediction

In practice, we cannot observe the whole predictor trajectory $X_i(t)$. Typically, it is assumed that $X_i(t)$ can only be realized at a discrete set of sampling points with additional measurement errors. Specifically, we observe the data $U_{il} = X_i(T_{il}) + \epsilon_{il}$, $T_{il} \in \mathcal{T}$; $l = 1, \ldots, N_i$, where $\epsilon_{il}$ are i.i.d. measurement errors with mean zero and finite variance $\sigma_u^2$, and $X_i(\cdot)$ are independent of $\epsilon_{il}$ for $i = 1, \ldots, n$. Depending on the number of observations $N_i$ within each curve, functional data are typically classified as sparse or dense; see Li and Hsing (2010) and Zhang and Wang (2016).

## 3 MODEL AND ESTIMATION

Many functions and parameters in the expressions given previously should be estimated from the data. We first use local linear smoothing to obtain the estimated mean function $\widehat{\mu}(t)$ and the estimated covariance function $\widehat{G}(s,t)$, which are well documented in the functional data analysis literature. For example, see Yao et al. (2005a,b), Li and Hsing (2010), and Zhang and Wang (2016) for the detailed calculations and we omit them here. Let the spectral decomposition of $\widehat{G}(s,t)$ be $\widehat{G}(s,t) = \sum_{j=1}^{\infty} \widehat{\kappa}_j \widehat{\phi}_j(s) \widehat{\phi}_j(t)$, where $\widehat{\kappa}_1 \geq \widehat{\kappa}_2 \geq \cdots \geq 0$ are the eigenvalues and $\{\widehat{\phi}_1, \widehat{\phi}_2, \ldots\}$ are corresponding eigenfunctions. Next, we use the principal components analysis through conditional expectation technique proposed by Yao et al. (2005a) to obtain an estimate of $\xi_{ij}$, denoted as $\widehat{\xi}_{ij}$. See Subsection S1.1 in the Supplementary Material for the details.

Let $\rho_\tau(e) = [\tau - \mathbf{1}\{e \leq 0\}]e$ be the check function, where $\mathbf{1}\{\cdot\}$ denotes the usual indicator function. The coefficients $a$ and $b_1, \ldots, b_{J_n}$ are estimated by

$$(\widehat{a}_{J_n}, \widehat{b}_{1,J_n}, \ldots, \widehat{b}_{J_n,J_n}) = \underset{a,b_1,\ldots,b_{J_n}}{\arg\min} \sum_{i=1}^{n} \rho_\tau \left( Y_i - a - \sum_{j=1}^{J_n} b_j \widehat{\xi}_{ij} \right).$$

Then, the resulting estimator of $b(t)$ is given by $\widehat{b}_{J_n}(t) = \sum_{j=1}^{J_n} \widehat{b}_{j,J_n} \widehat{\phi}_j(t)$.

Let $(Y_0, X_0)$ be an independent copy of $\{(Y_i, X_i)\}_{i=1}^n$. Similar to $X_i(t)$, $X_0(t)$ is observed with measurement error, i.e., only $\{(T_{0l}, U_{0l})\}_{l=1}^{N_0}$ with $i = 0$ is observed. As a result, the $\tau$th conditional quantile of $Y_0$ given $X_0$ is predicted by a plug-in method:

$$\widehat{Q}_{\tau,J_n}(X_0) = \widehat{a}_{J_n} + \int_{\mathcal{T}} \widehat{b}_{J_n}(t) \left\{ \widehat{X}_{0,J_n}(t) - \widehat{\mu}(t) \right\} dt = \widehat{a}_{J_n} + \sum_{j=1}^{J_n} \widehat{b}_{j,J_n} \widehat{\xi}_{0j}.$$

Here, $\widehat{\xi}_{0j}$ is obtained similarly to $\widehat{\xi}_{ij}$, and $\widehat{X}_{0,J_n}(t) = \widehat{\mu}(t) + \sum_{j=1}^{J_n} \widehat{\xi}_{0j} \widehat{\phi}_j(t)$. The predicted conditional quantile $\widehat{Q}_{\tau,J_n}(X_0)$ depends on the number of principal components $J_n$, and thus the prediction performance varies with $J_n$. In the case where $\{X_i(t)\}_{i=0}^n$ are observed at dense discrete points without measurement errors, Kato (2012) established the minimax optimal rates of convergence of $\widehat{b}_{J_n}(t)$ and $\widehat{Q}_{\tau,J_n}(X_0)$ with a proper choice of $J_n$.

## 4. Model Averaging for FLQR Model

### 4.1. Weighted Average Quantile Prediction

The choice of $J_n$ is usually made based on a model selection criterion such as FVE, AIC, or BIC. Let the prediction of $Q_\tau(X_0)$ from a fixed $J_n$ choice be $\widehat{Q}_{\tau,J_n}(X_0)$, where $J_n \in \mathcal{J}$ for a candidate set $\mathcal{J}$. Typically, $\mathcal{J} = \{J_L, J_L + 1, \ldots, J_U\}$, where $J_L$ and $J_U$ are lower and upper bounds, respectively. Let $w_{J_n}$ be a weight assigned to the model with $J_n$ eigenbases. Let $\mathbf{w}$ be a vector formed by all such weights $w_{J_n}$. For example, $\mathbf{w} = (w_{J_L}, w_{J_L+1}, \ldots, w_{J_U})^\top$ if $\mathcal{J} = \{J_L, J_L + 1, \ldots, J_U\}$. Let $\mathcal{W} = \{\mathbf{w}: w_{J_n} \geq 0, J_n \in \mathcal{J}, \sum_{J_n \in \mathcal{J}} w_{J_n} = 1\}$. The model averaging prediction with weights $\mathbf{w}$ is $\widehat{Q}_\tau(X_0, \mathbf{w}) = \sum_{J_n \in \mathcal{J}} w_{J_n} \widehat{Q}_{\tau,J_n}(X_0)$.

We define the final prediction error (FPE, or the out-of-sample prediction error) used by Lu and Su (2015) and Wang et al. (2023) as $\text{FPE}_n(\mathbf{w}) = E\left\{\rho_\tau\left(Y_0 - \widehat{Q}_\tau(X_0, \mathbf{w})\right)\Big|\mathcal{D}_n\right\}$, where $\mathcal{D}_n = \{(Y_i, T_{il}, U_{il}): l = 1, \ldots, N_i; i = 1, \ldots, n\}$ is the observed sample. Let $F(\cdot|X_i)$ denote the conditional distribution function of $\varepsilon_i$ given $X_i$. Using the identity (Knight 1998)

$$\rho_\tau(u - v) - \rho_\tau(u) = -v\psi_\tau(u) + \int_0^v [\mathbf{1}\{u \leq s\} - \mathbf{1}\{u \leq 0\}]\, ds, \tag{4.4}$$

where $\psi_\tau(u) = \tau - \mathbf{1}\{u \leq 0\}$, we have $\text{FPE}_n(\mathbf{w}) - E\{\rho_\tau(\varepsilon_0)\} \geq 0$; see the detailed discussions in the Supplementary Material.

Since $E\{\rho_\tau(\varepsilon_0)\}$ is unrelated to $\mathbf{w}$, minimizing $\text{FPE}_n(\mathbf{w})$ is equivalent to minimizing the following excess final prediction error (EFPE):

$$\text{EFPE}_n(\mathbf{w}) = \text{FPE}_n(\mathbf{w}) - E\{\rho_\tau(\varepsilon_0)\}.$$

From Equation (S1.2), it is easy to see that $\text{EFPE}_n(\mathbf{w}) \geq 0$ for each $\mathbf{w} \in \mathcal{W}$. A similar predictive risk of the quantile regression model is considered by Giessing and He (2019). In addition, EFPE is closely related to the usual mean squared prediction error $R_n(\mathbf{w}) = E[\{\widehat{Q}_\tau(X_0, \mathbf{w}) - Q_\tau(X_0)\}^2|\mathcal{D}_n]$.

To see it, under Conditions 2 and 5 in the Supplementary Material, it is easy to prove that $\underline{C} \leq \text{EFPE}_n(\mathbf{w})/R_n(\mathbf{w}) \leq \bar{C}$ for each $\mathbf{w} \in \mathcal{W}$, where $0 < \underline{C} \leq \bar{C} < \infty$ are defined in Condition 5.

## 4.2. Weight Choice Criterion

We use $K$-fold cross-validation to choose the weights. Specifically, we divide the dataset into $K \geq 2$ groups such that the sample size of each group is $M = \lfloor n/K \rfloor$, where $\lfloor \cdot \rfloor$ denotes the greatest integer less than or equal to $\cdot$. When $K = n$, we have $M = 1$, and $K$-fold cross-validation becomes leave-one-curve-out cross-validation used by Lu and Su (2015) and Wang et al. (2023), which is not computationally feasible when $n$ is too large. For each fixed $J_n \in \mathcal{J}$, we consider the $k$th step of the $K$-fold cross-validation ($k = 1, \ldots, K$), where we leave out the $k$th group and use all of the remaining observations to obtain the intercept and the slope function estimates $\widehat{a}_{J_n}^{[-k]}$ and $\widehat{b}_{J_n}^{[-k]}(t)$, respectively. Now, for $i = (k-1)M+1, \ldots, kM$, we estimate $Q_\tau(X_i)$ by

$$\widehat{Q}_{\tau,J_n}^{[-k]}(X_i) = \widehat{a}_{J_n}^{[-k]} + \int_{\mathcal{T}} \widehat{b}_{J_n}^{[-k]}(t) \left\{ \widehat{X}_{i,J_n}(t) - \widehat{\mu}(t) \right\} dt,$$

where $\widehat{X}_{i,J_n}(t) = \widehat{\mu}(t) + \sum_{j=1}^{J_n} \widehat{\xi}_{ij} \widehat{\phi}_j(t)$. Note that $\widehat{\mu}(t)$ and $\widehat{X}_{i,J_n}(t)$ are obtained by using all data. Then, the weighted average prediction with weights $\mathbf{w}$ is $\widehat{Q}_\tau^{[-k]}(X_i, \mathbf{w}) = \sum_{J_n \in \mathcal{J}} w_{J_n} \widehat{Q}_{\tau,J_n}^{[-k]}(X_i)$. The $K$-fold cross-validation criterion is formulated as $\text{CV}_K(\mathbf{w}) = \frac{1}{n} \sum_{k=1}^{K} \sum_{m=1}^{M} \rho_\tau \left( Y_{(k-1)M+m} - \widehat{Q}_\tau^{[-k]}(X_{(k-1)M+m}, \mathbf{w}) \right)$. The resulting weight vector is obtained as

$$\widehat{\mathbf{w}} = \arg\min_{\mathbf{w} \in \mathcal{W}} \text{CV}_K(\mathbf{w}). \tag{4.5}$$

Consequently, the proposed model average $\tau$th conditional quantile prediction of $Y_0$ given $X_0$ is $\widehat{Q}_\tau(X_0, \widehat{\mathbf{w}})$. The constrained minimization problem (4.5) can be reformulated as a linear programming problem provided in Subsection S1.3 of the Supplementary Material. The MATLAB code for our method is available from `https://github.com/wcxstat/MAFLQR`.

### 4.3. Choice of the Candidate Set

In practice, the candidate set $\mathcal{J}$ can be chosen heuristically, coupled with a model selection criterion. Let $d$ be a fixed small positive integer, and consider the set $\mathcal{J} = \{j \in \mathbb{N}: |\widehat{J}_n - j| \leq d\}$, where $\widehat{J}_n$ can be selected by the FVE, AIC, or BIC criterion. Specifically, following Chen and Müller (2012), FVE is defined as $\text{FVE}_{J_n} = \sum_{j=1}^{J_n} \widehat{\kappa}_j / \sum_{j=1}^{\infty} \widehat{\kappa}_j$, and following Kato (2012) and Lu and Su (2015), for model $J_n$, the AIC and BIC are respectively defined as $\text{AIC}_{J_n} = 2n \log \left\{ \frac{1}{n} \sum_{i=1}^{n} \rho_\tau \left( Y_i - a - \sum_{j=1}^{J_n} \widehat{b}_{j,J_n} \widehat{\xi}_{ij} \right) \right\} + 2(J_n + 1)$ and $\text{BIC}_{J_n} = 2n \log \left\{ \frac{1}{n} \sum_{i=1}^{n} \rho_\tau \left( Y_i - a - \sum_{j=1}^{J_n} \widehat{b}_{j,J_n} \widehat{\xi}_{ij} \right) \right\} + (J_n + 1) \log n$. Given a threshold $\gamma$ (e.g., 0.90 or 0.95), FVE selects $\widehat{J}_n$ as the smallest integer that satisfies $\text{FVE}_{J_n} \geq \gamma$; AIC and BIC select $\widehat{J}_n$ by minimizing $\text{AIC}_{J_n}$ and $\text{BIC}_{J_n}$, respectively.

We refer to the resultant procedure as the *d*-divergence model averaging method. When $d = 0$, it reduces to the selection criterion without model averaging. In Section 6, we compare the numerical performance of various *d*-divergence model averaging methods under different settings of $\mathcal{J}$. For the following theoretical studies, we assume $\mathcal{J}$ in our method is predetermined.

## 5. Asymptotic Results

### 5.1. Asymptotic Weight Choice Optimality

In this subsection, we present an important result on the asymptotic optimality of the selected weights. Denote by $|\mathcal{J}|$ the cardinality of the set $\mathcal{J}$, which is also called the number of candidate models. All limiting processes are studied with respect to $n \to \infty$.

**Theorem 1.** *Suppose that Conditions 1–4 in the Supplementary Material hold. Then,*

$$\frac{\text{EFPE}_n(\widehat{\mathbf{w}})}{\inf_{\mathbf{w} \in \mathcal{W}} \text{EFPE}_n(\mathbf{w})} = 1 + O_p \left\{ \eta_n^{-1} n \left( c_n + g_n + n^{-1/2} |\mathcal{J}|^{1/2} \right) \right\} = 1 + o_p(1).$$

Theorem 1 shows that the selected weight vector $\widehat{\mathbf{w}}$ is asymptotically optimal in the sense that its EFPE is asymptotically identical to that of the infeasible best weight vector to minimize $\text{EFPE}_n(\mathbf{w})$.

# 5 ASYMPTOTIC RESULTS

The asymptotic optimality in Theorem 1 is better than the existing results in the literature such as Hansen and Racine (2012), Lu and Su (2015), and Zhang, Chiou, and Ma (2018) because the former provides a convergence rate but the latter does not. The rate is determined by $\eta_n$, the convergence rate $c_n$ of estimators of each candidate model, the difference $g_n$ between the regular prediction $\widehat{Q}_\tau(X_{(k-1)M+m}, \mathbf{w})$ and the leave-$M$-out prediction $\widehat{Q}_\tau^{[-k]}(X_{(k-1)M+m}, \mathbf{w})$, and the number of candidate models $|\mathcal{J}|$.

Note that Lu and Su (2015) and Wang et al. (2023) established asymptotic optimality in terms of minimizing $\mathrm{FPE}_n(\mathbf{w})$ instead of $\mathrm{EFPE}_n(\mathbf{w})$. We present a similar result in the following corollary.

**Corollary 1.** *Suppose that Conditions 1–3 hold. If $n^{-1}|\mathcal{J}| \to 0$ as $n \to \infty$, then*

$$\frac{\mathrm{FPE}_n(\widehat{\mathbf{w}})}{\inf_{\mathbf{w} \in \mathcal{W}} \mathrm{FPE}_n(\mathbf{w})} = 1 + O_p\left(c_n + g_n + n^{-1/2}|\mathcal{J}|^{1/2}\right) = 1 + o_p(1).$$

Corollary 1 can be proven using analogous arguments as those used in the proof of Theorem 1, and thus we omit its proof. Since $\inf_{\mathbf{w} \in \mathcal{W}} \mathrm{FPE}_n(\mathbf{w}) \geq E\{\rho_\tau(\varepsilon_0)\} > 0$ from Equation (S1.2), the asymptotic optimality in Corollary 1 is implied by that in Theorem 1 with $n^{-1}|\mathcal{J}| \to 0$ as $n \to \infty$. Therefore, the asymptotic optimality in Theorem 1 improves the results found in the studies of Lu and Su (2015) and Wang et al. (2023). This is an important contribution to model averaging in quantile regression.

The results of Theorem 1 and Corollary 1 hold for a single chosen quantile $\tau$. Suppose we have a given subset $\mathcal{U}$ of $(0, 1)$ that is away from 0 and 1. Examples of such subsets include $\mathcal{U}$ are $\mathcal{U} = \{\tau_1, \ldots, \tau_d\}$ with $0 < \tau_1 < \cdots < \tau_d < 1$ and $\mathcal{U} = [\tau_L, \tau_U]$ with $0 < \tau_L < \tau_U < 1$. If Conditions 1–3 hold uniformly for $\tau \in \mathcal{U}$, it is easy to show that the asymptotic optimality in Theorem 1 and Corollary 1 hold uniformly for $\tau \in \mathcal{U}$.

We finalize this subsection by discussing the advantage of using EFPE over FPE. Observe that for $J_n \in \mathcal{J}$, $E\left\{\left|\widehat{Q}_{\tau,J_n}(X_0) - Q_\tau(X_0)\right|\Big|\mathcal{D}_n\right\} \leq |\widehat{a}_{J_n} - a| + \int_{\mathcal{T}} \left|\widehat{b}_{J_n}(t) - b(t)\right| E\left\{\left|\widehat{X}_{0,J_n}(t) - \widehat{\mu}(t)\right|\Big|\mathcal{D}_n\right\} dt +$

$\int_{\mathcal{T}} |b(t)| \left[ E\left\{ \left\| \widehat{X}_{0,J_n}(t) - X_0(t) \right\| \big| \mathcal{D}_n \right\} + |\widehat{\mu}(t) - \mu(t)| \right] dt$. Let $J^*$ correspond to a true model (we allow $J^* = \infty$), i.e., $b(t) = \sum_{j=1}^{J^*} b_j \phi_j(t)$. Note that this may not be the only true model, and any model containing model $J^*$ is also true. We consider two scenarios: (i) There is at least one true model (say $J_n$) in the set of candidate models. (ii) All candidate models are misspecified, but one of the candidate models (say $J_n$) converges to a true model with a rate.

Under the scenario (i) or (ii), it can be proven that $\widehat{a}_{J_n} - a = o_p(1)$, $\widehat{b}_{J_n}(t) - b(t) = o_p(1)$, $E\{|\widehat{X}_{0,J_n}(t) - X_{0,J_n}(t)| | \mathcal{D}_n\} = o_p(1)$, and $\widehat{\mu}(t) - \mu(t) = o_p(1)$ hold uniformly for $t \in \mathcal{T}$ under certain conditions, where $X_{0,J_n}(t) = \mu(t) + \sum_{j=1}^{J_n} \xi_{0j} \phi_j(t)$; see, e.g., Yao et al. (2005a), Kato (2012), and Li et al. (2022). This result implies that $E\{|\widehat{Q}_{\tau,J_n}(X_0) - Q_\tau(X_0)| | \mathcal{D}_n\} \to 0$ in probability from (5.1), i.e., $Q_\tau(X_0)$ can be consistently predicted by a single candidate model. This, along with Equation (S1.2), implies that under scenario (i) or (ii), $\inf_{\mathbf{w} \in \mathcal{W}} \text{EFPE}_n(\mathbf{w}) \to 0$ in probability. Compared to $\text{EFPE}_n(\mathbf{w})$, $\text{FPE}_n(\mathbf{w})$ contains an additional term $E\{\rho_\tau(\varepsilon_0)\}$, which is a positive number unrelated to $n$ or $\mathbf{w}$. Therefore, under scenario (i) or (ii), $E\{\rho_\tau(\varepsilon_0)\}$ is the dominant term in FPE, which makes the asymptotic optimality built based on FPE no sense. This phenomenon has not been found in existing literature such as Lu and Su (2015) and Wang et al. (2023).

## 5.2. Estimation Consistency

Define the model averaging estimators of $a$ and $b(t)$ with weights $\mathbf{w}$ as

$$\widehat{a}_{\mathbf{w}} = \sum_{J_n \in \mathcal{J}} w_{J_n} \widehat{a}_{J_n} \quad \text{and} \quad \widehat{b}_{\mathbf{w}}(t) = \sum_{J_n \in \mathcal{J}} w_{J_n} \widehat{b}_{J_n}(t),$$

respectively. Let $J^* < \infty$ correspond to a true model. Assume there exists a series $d_n \to 0$ such that

$$\widehat{a}_{J^*} - a = O_p(d_n) \quad \text{and} \quad \int_{\mathcal{T}} \left\{ \widehat{b}_{J^*}(t) - b(t) \right\}^2 dt = O_p(d_n^2). \tag{5.6}$$

This states that the estimators $\widehat{a}_{J^*}$ and $\widehat{b}_{J^*}(t)$ are consistent, which are verified in the literature (e.g., see Li et al. 2022, Kato 2012). In the ideal case that the functional covariate $X_i(t)$ is fully observed, $d_n = n^{-1/2}$ for the fixed $J^*$. In more common cases where the curves are observed at some discrete points, $d_n$ may be slower than $n^{-1/2}$ since $\mu(t)$ and $G(s,t)$ should be estimated first by some nonparametric smoothing methods.

We now describe the performance of the weighted estimation results when there is at least one true model among the candidate models.

**Theorem 2.** *When there is at least one true model among candidate models, under Conditions 1–3 and 5–7 in the Supplementary Material, and the assumption $n^{-1}|\mathcal{J}| \to 0$ as $n \to \infty$, the weighted average estimators $\widehat{a}_{\widehat{\mathbf{w}}}$ and $\widehat{b}_{\widehat{\mathbf{w}}}(t)$ satisfy*

$$(\widehat{a}_{\widehat{\mathbf{w}}} - a)^2 + \int_{\mathcal{T}} \left\{\widehat{b}_{\widehat{\mathbf{w}}}(t) - b(t)\right\}^2 dt = O_p\left(d_n^2 + c_n + g_n + n^{-1/2}|\mathcal{J}|^{1/2}\right),$$

*where $c_n$ and $g_n$ are defined in Conditions 1 and 3, respectively.*

Theorem 2 provides a convergence rate for $\widehat{a}_{\widehat{\mathbf{w}}}$ and $\widehat{b}_{\widehat{\mathbf{w}}}(t)$ when at least one true model is included in the candidate models. In this situation, this rate consists of two parts: the convergence rate $d_n$ of $(\widehat{a}_{J^*}, \widehat{b}_{J^*})$ under the true model, and an additional term $c_n + g_n + n^{-1/2}|\mathcal{J}|^{1/2}$ that arises from the estimated weights. However, this rate may not be optimal and could potentially be improved with further efforts. This is an interesting direction for future research. It should be noted that if (5.6) and Conditions 1–3 and 7 hold uniformly for $\tau \in \mathcal{U}$, we can also show that the result in Theorem 2 holds uniformly for $\tau \in \mathcal{U}$.

## 6. Simulation Studies

In this section, we conduct simulation studies to evaluate the finite sample performance of the proposed model averaging quantile prediction and the estimation of the slope function. We consider two simulation designs. In the first design, all candidate models are misspecified, whereas in the second, there are true models included in the set of candidate models. To save space, the simulation studies for the second design are provided in Subsection S4.2 of the Supplementary Material.

### 6.1. Comparison Methods

We compare the proposed procedure with several model selection methods and two existing model averaging methods, namely the smoothed AIC (SAIC) and smoothed BIC (SBIC) proposed by Buckland et al. (1997). The model selection methods include FVE with $\gamma \in \{0.90, 0.95\}$, AIC, and BIC, as described in Subsection 4.3. The SAIC and SBIC model averaging methods assign weights

$$w_{\text{AIC}, J_n} = \frac{\exp\left(-\frac{1}{2}\text{AIC}_{J_n}\right)}{\sum_{J_n \in \mathcal{J}} \exp\left(-\frac{1}{2}\text{AIC}_{J_n}\right)} \quad \text{and} \quad w_{\text{BIC}, J_n} = \frac{\exp\left(-\frac{1}{2}\text{BIC}_{J_n}\right)}{\sum_{J_n \in \mathcal{J}} \exp\left(-\frac{1}{2}\text{BIC}_{J_n}\right)}$$

to model $J_n$, respectively, where $\text{AIC}_{J_n}$ and $\text{BIC}_{J_n}$ are defined in Subsection 4.3.

### 6.2. Simulation Design I

Set $n_T = n + n_0$ as the sample size of a dataset, among which the first $n$ dataset are used as training data and the other $n_0 = 100$ as test data. Consider sparse designs for $X_i$. The observations $\{(T_{il}, U_{il}, X_i)\}$ are sampled by the steps in Subsection S4.1 of the Supplementary Material. Then, we generate the response observation $Y_i$ by the following heteroscedastic functional linear model

$$Y_i = \theta \int_0^1 b^{[1]}(t) X_i(t)\, dt + \sigma(X_i) e_i, \tag{6.7}$$

where $\sigma(X_i) = a^{[2]} + \int_0^1 b^{[2]}(t) X_i(t)\, dt$, $e_i$ are sampled from $N(0, 1)$, and $b^{[1]}$, $a^{[2]}$, and $b^{[2]}$ are given in Subsection S4.1 of the Supplementary Material.

It is easy to see that (6.7) leads to an FLQR model of the form (3.1) with $a = F_e^{-1}(\tau) a^{[2]}$ and $b(t) = \theta b^{[1]}(t) + F_e^{-1}(\tau) b^{[2]}(t)$, where $F_e^{-1}(\cdot)$ is the quantile function of the distribution of $e_i$. As in Lu and Su (2015), we define the population $R^2$ as $R^2 = [\text{Var}(Y_i) - \text{Var}\{\sigma(X_i) e_i\}] / \text{Var}(Y_i)$. We consider $\tau = 0.5, 0.05$ and different choices of $\theta$ such that $R^2$ ranges from 0.1 to 0.9 with an increment of 0.1. To evaluate each method, we compute the excess final prediction error. We do this by computing averages across 200 replications.

We begin by comparing the performance of the $d$-divergence model averaging under different settings of $\mathcal{J}$ described in Subsection 4.3. Specifically, we let $d = 0, 1, 2, 4$, and $\widehat{J}_n$ be selected by FVE with $\gamma \in \{0.90, 0.95\}$, AIC, or BIC. We fix $n = 300$ and $R^2 = 0.5$. For each design, we present the results with $K = 2, 4$ for the $K$-fold cross-validation method to obtain the optimal weights. Figure A.1 shows the boxplots of EFPE with $\tau = 0.05$. The results for $\tau = 0.5$ are displayed in Figure S.9 of the Supplementary Material. When $\tau = 0.5$, compared to each model selection method, the corresponding model averaging method provides little further improvement. However, when $\tau = 0.05$, the corresponding model averaging method can indeed further reduce the excess final prediction error compared to each model selection method.

Next, we compare the proposed model averaging prediction with the other six model averaging and selection methods. We choose $d$-divergence model averaging with $d = 4$ coupled with the FVE(0.90) criterion as our model averaging method. We consider $n = 100, 200, 400$ and $\tau = 0.5, 0.05$, and $K$ is fixed to be 4. In each simulation setting, we normalize the EFPEs of the other six methods by a division of the EFPE of our method. The results are presented in Figure A.2. The performance of different methods becomes more similar when $R^2$ is large. When $\tau = 0.5$, no method clearly dominates the others. Our method is the best in most cases for $n = 100, 200$. The performances of AIC and BIC seem to be the worst. When $\tau = 0.05$, it is clear that our method significantly dominates all the other six methods, and SBIC and FVE(0.90) seem to be the second and third best, respectively when $R^2 \leq 0.5$. The reason why our method has more advantages than other methods for extreme quantile cases may be that the data at extreme quantiles are usually sparse, the prediction for extreme quantiles is more challenging, and model averaging can address this challenge by using more information from more models. To further show this performance comprehensively, we conduct additional simulation studies with $\tau = 0.2$ and $0.1$. The results are displayed in Figure S.10 of the Supplementary Material. From Figures A.2 and S.10, the advantage of model averaging over the other six methods increases as $\tau$ decreases. Overall, our study illustrates the advantage of model averaging over model selection in terms of prediction for extreme quantiles.

## 7. Experiments

Accurate quantile forecasting is critical for financial risk management, including VaR estimation and capital allocation decisions (Gao et al. 2019, Jin et al. 2021, Lai et al. 2024, Liu et al. 2022). The extreme volatility, heavy-tailed return distributions, and unique dynamics of cryptocurrencies make VaR forecasting a highly complex and challenging task. This study provides a conditional quantile model averaging prediction framework leveraging functional data to address this issue. This section will examine the predictive performance of the proposed MA with different methods in Section 6.1, including SAIC, SBIC, FVE90, FVE95, AIC and BIC.

This study examines the price data from ten major cryptocurrencies: Cardano (ADA), Avalanche (AVAX), Bitcoin Cash (BCH), Bitcoin (BTC), Ethereum (ETH), Chainlink (LINK), Litecoin (LTC), Uniswap (UNI), Stellar (XLM), and Ripple (XRP). The data set was collected from June 1, 2022 and May 31, 2024 on the exchange platform Bitstamp. We use our method to, based on hourly log-returns of the price of Bitcoin on a given day, predict the $\tau$th quantile of the minimum hourly log-return of Bitcoin price the next day. To reduce the temporal dependence existing in the data, as in Girard et al. (2022), we construct our sample of data by keeping a gap of one day between observations. Specifically, in the $i$th sample $(Y_i, X_i)$, the functional covariate $X_i$ is the curve of hourly log-returns on day $2i - 1$ and the scalar response $Y_i$ is the minimum hourly log-return on

https://www.cryptodatadownload.com/data/

day 2$i$. The left panels of Figures A.3-A.4 provide hourly log-returns of Bitcoin price, and the right panels plot the histogram of minimum hourly log-return of Bitcoin price. Obviously the distribution of minimum hourly log-return is highly right skewed; hence a simple statistic like sample mean cannot adequately describe the (conditional) distribution of $Y_i$ given $X_i$. We use the sample during June 1 2022 to May 31 2023 as the training data to build the FLQR model and use the remaining samples as the test data to evaluate the prediction performance. To better assess prediction accuracy, we repeat this $B = 200$ times based on random partitions of the data set.

To evaluate each method, we calculate the final prediction error, which are calculated as

$$\text{FPE}(\mathbf{w}) = \frac{1}{\lfloor 0.3n \rfloor B} \sum_{r=1}^{B} \sum_{i=1}^{\lfloor 0.3n \rfloor} \rho_\tau \left( Y_{0,i}^{(r)} - \widehat{Q}_\tau(X_{0,i}^{(r)}, \mathbf{w})^{(r)} \right),$$

where $\{(Y_{0,i}^{(r)}, X_{0,i}^{(r)})\}_{i=1}^{\lfloor 0.3n \rfloor}$ is the test data and $\widehat{Q}_\tau(\cdot, \mathbf{w})^{(r)}$ denotes the prediction based on the training data in the $r$th partition. We consider the $\tau$th quantile prediction with $\tau = 0.05$ and $0.01$. By Equation (S1.2), the above FPE also measures the EFPE. To perform prediction using model averaging, we choose the $d$-divergence model averaging described in Subsection 4.3, where $d = 0, 2, 4, 8, 10, 16, 22$ and $\widehat{J}_n$ is selected by FVE(0.90), FVE(0.95), AIC, and BIC.

Figures A.5-A.6 illustrates the final prediction errors for $K = 2$ under ten cryptocurrencies. Note that the result of $d = 0$ corresponds to that of model selection with $\widehat{J}_n$. Figures A.5-A.6 suggest that for each $\tau$, the predictive ability based on selecting the number of principal components using model selection criteria is unstable across different assets. For some assets, BIC performs well because it selects a small number of principal components, while for others, FVE90 performs better. ACI consistently performs the worst, overestimating the numbers of principal components. However, these results are much improved when using model averaging with $d = 10$ or larger. The data example demonstrates that model averaging improves the prediction performance and can protect against potential prediction loss caused by a single selection criterion.

Besides, we compare our method with the other six model averaging and selection methods described in Section 6. We fix $K = 2$ and choose $\mathcal{J}$ in Subsection 4.3 with $d = 8$ and $\widehat{J}_n$ selected by BIC as the candidate set. Figures A.7–A.8 report the boxplots of FPEs based on the seven methods in excess return forecasting. It is observed that the proposed QMA is consistently ranked first among other competing methods. For instance, for $\tau = 0.01$ and BCH, our method outperforms SAIC and SBIC, achieving nearly a twofold improvement. These improvements may come from our asymptotically optimal weight selection criterion. Besides, the sAIC and sBIC methods outperform the related model selections AIC and BIC in most cases. This highlights the merits of reducing model uncertainty and increasing the robustness of model averaging. Overall, these empirical results highlight the adavantages of our method, especially for extreme quantiles.

Furthermore, we evaluate these VaR forecasts through calibration tests, which are typically employed for assessing the accuracy of quantile predictions. For a probability level $\tau$, a quantile forecast $\widehat{Q}_t$ is said to be unconditionally calibrated if $\text{Hit}_t = \tau - I(y_t \leq \widehat{Q}_t)$ has zero unconditional expectation, and conditionally calibrated if $\text{Hit}_t$ has zero conditional expectation. Following Taylor (2020), we assess unconditional calibration using a binomial test that checks if the mean of $\text{Hit}_t$ significantly differs from zero. Additionally, we consider other classical VaR backtesting methods, including the proportion of failures test (POF), conditional coverage independence test (CC), and time between failures mixed test (TBF). Table A.1 summarizes the test results for the ten cryptocurrencies at the 1% and 5% probability levels. The reported values correspond to the number of indices for which the test was not rejected at the given significance level, with larger values indicating better performance. Table A.1 demonstrates that the proposed MA approach outperforms other competing methods, particularly at the 5% level.

## 8. Conclusion

We develop a *K*-fold cross-validation model averaging for the FLQR model. When there is no true model included in candidate models, the weight vector chosen by our method is asymptotically optimal in terms of minimizing the excess final prediction error, whereas when there is at least one true model included, the model averaging parameter estimates are consistent. Our simulations and empirical applications to tail risks forecasting in cryptocurrency markets indicate that the proposed method outperforms the other model averaging and selection methods.

There are some open questions for future research. First, the asymptotic behavior of the selected weights and post-model-averaging inference for quantile regression were not studied in this work. While post-model-averaging inference for conditional means has been explored in Zhang and Liu (2019) and Yu et al. (2024), the quantile case is unknown. Additionally, extending model averaging prediction to multimodal data beyond functional data is an interesting direction.

# Appendix

## A1. Simulation Results

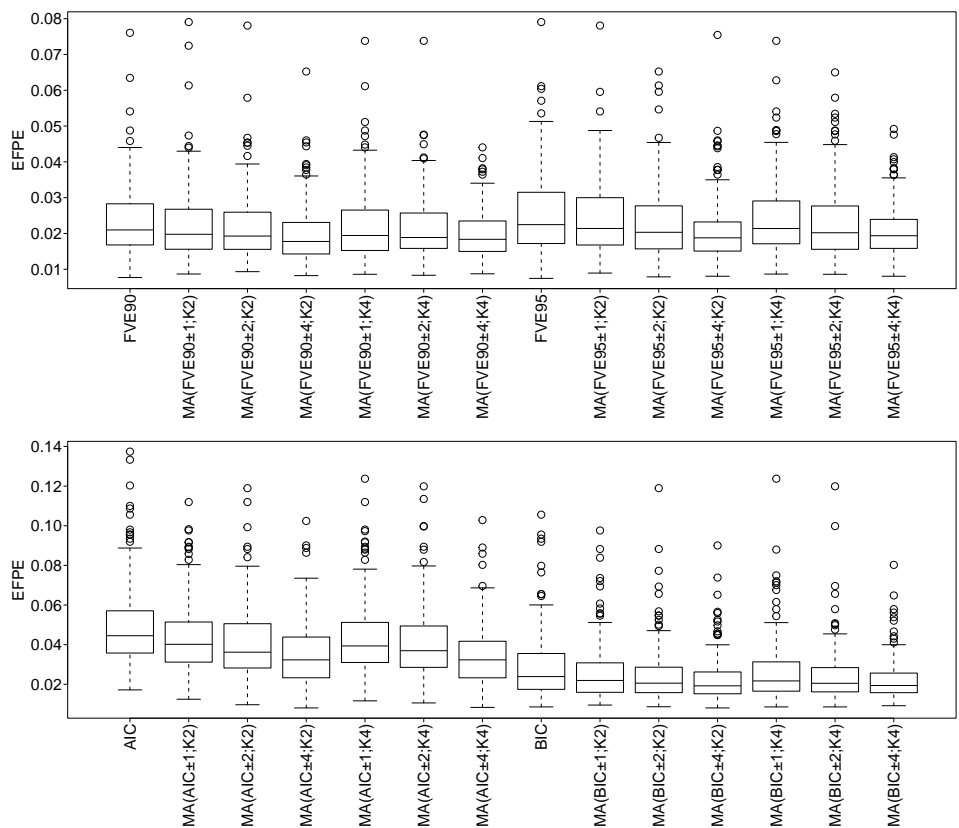

**Figure A.1** **Boxplots of EFPE in simulation design I with** $\tau = 0.05$**. FVE90, FVE with** $\gamma = 0.90$**; FVE95, FVE with** $\gamma = 0.95$**; MA(FVE90**$\pm\alpha$**,** $K\beta$**), model averaging method with** $\widehat{J}_n$ **determined by FVE90, with** $d = \alpha$**, and with weights selected by** $\mathrm{CV}_{K=\beta}$**; MA(FVE95**$\pm\alpha$**,** $K\beta$**), MA(AIC**$\pm\alpha$**,** $K\beta$**), and MA(BIC**$\pm\alpha$**,** $K\beta$**) have similar definitions.**

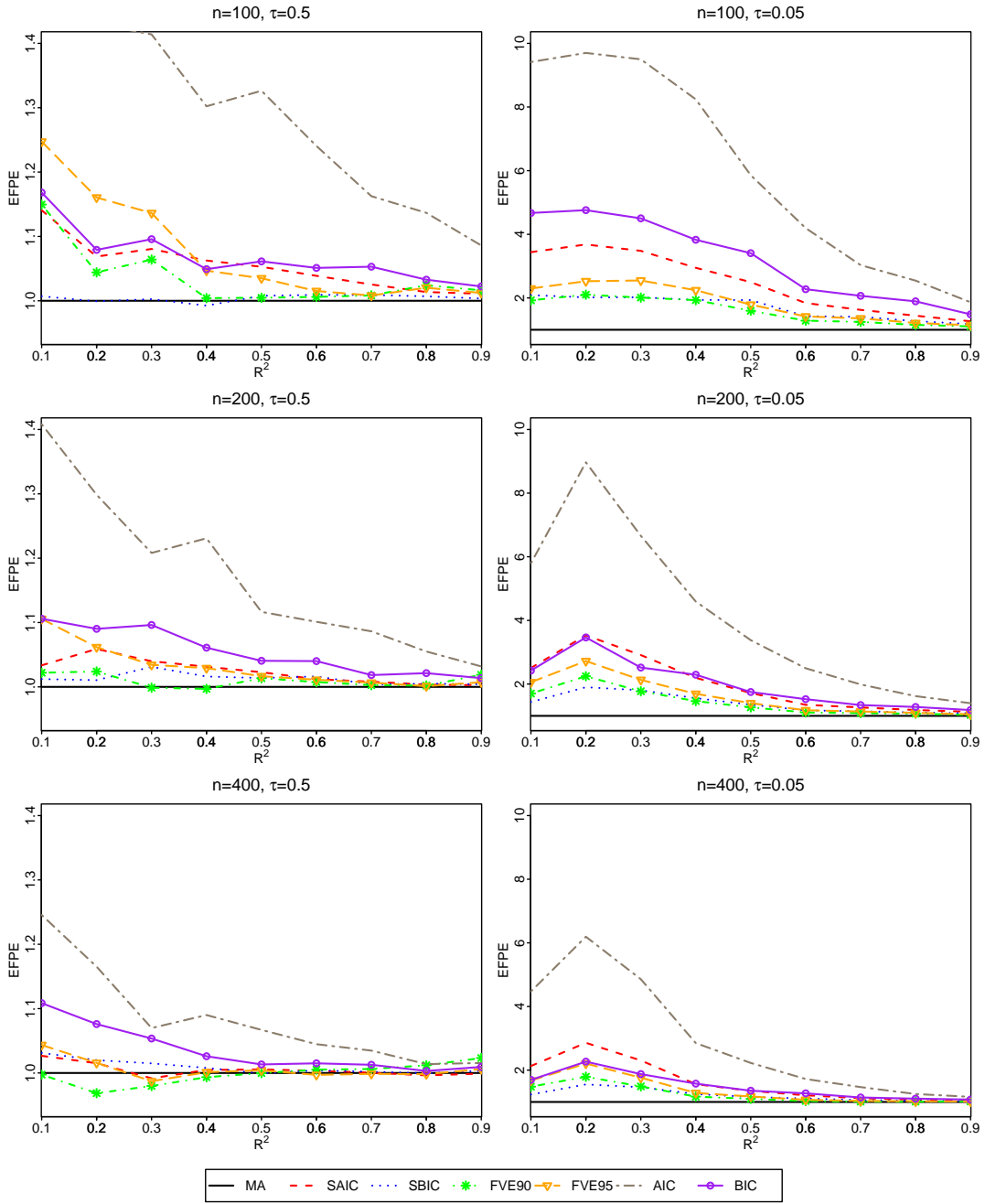

**Figure A.2**   Normalized EFPE in simulation design I for $\tau = 0.5$ and $0.05$. **MA denotes MA(FVE90$\pm$4, $K$4).**

## A2.   Empirical Results

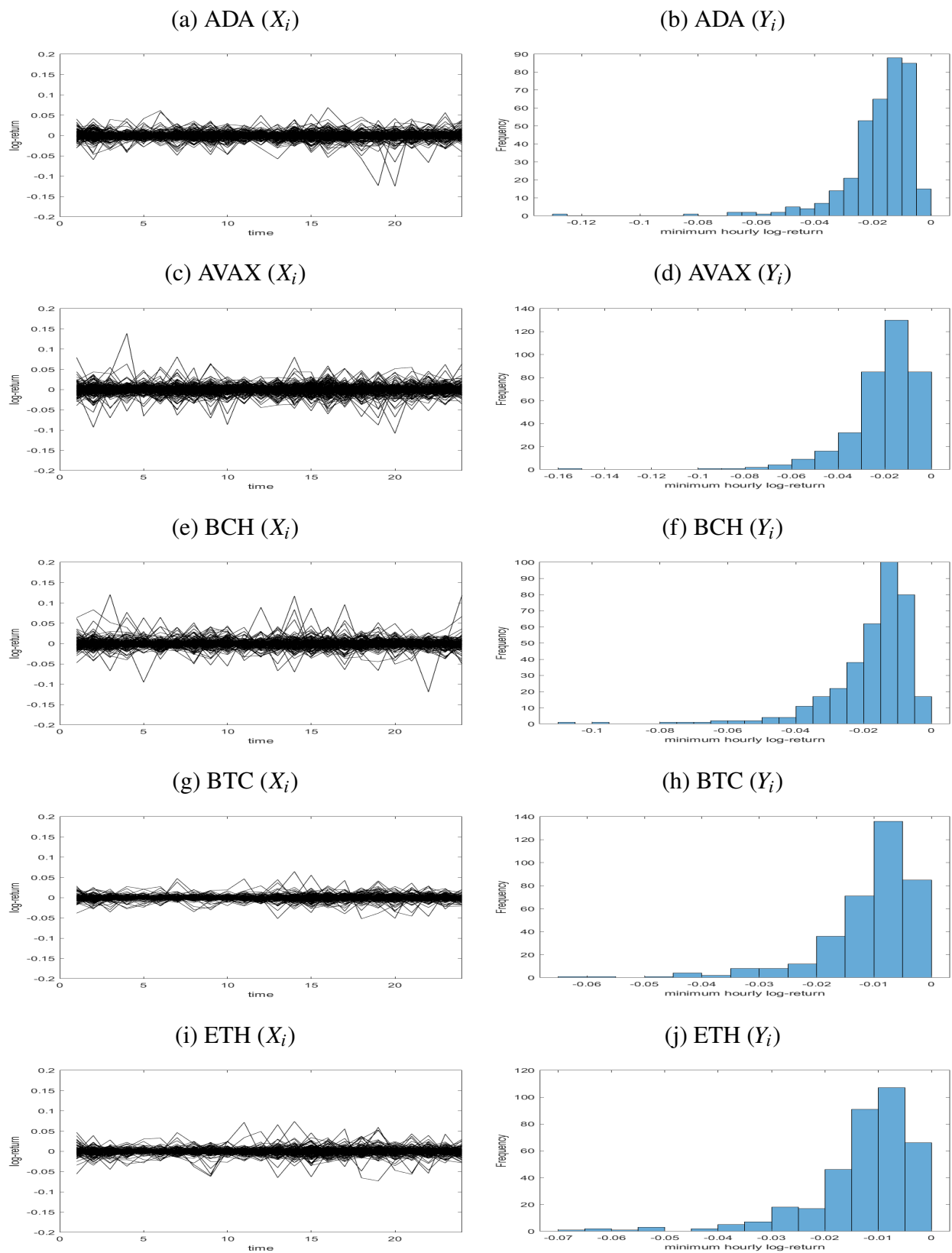

**Figure A.3** Cryptocurrency returns. Left: Hourly log-returns of cryptocurrency prices. Right: Histogram of the minimum hourly log-return of cryptocurrency prices.

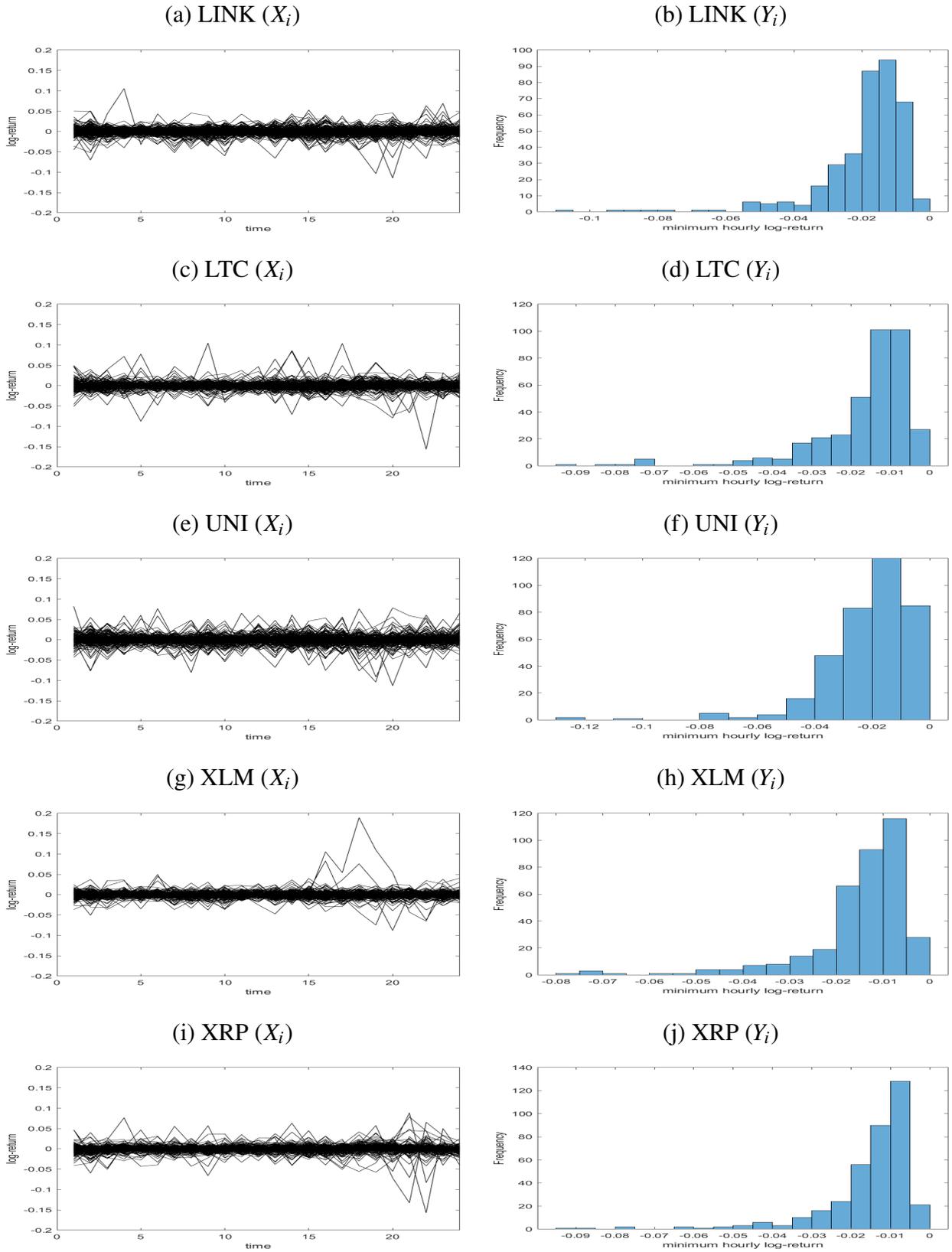

**Figure A.4** Cryptocurrency returns. Left: Hourly log-returns of cryptocurrency prices. Right: Histogram of the minimum hourly log-return of cryptocurrency prices.

## A2 EMPIRICAL RESULTS

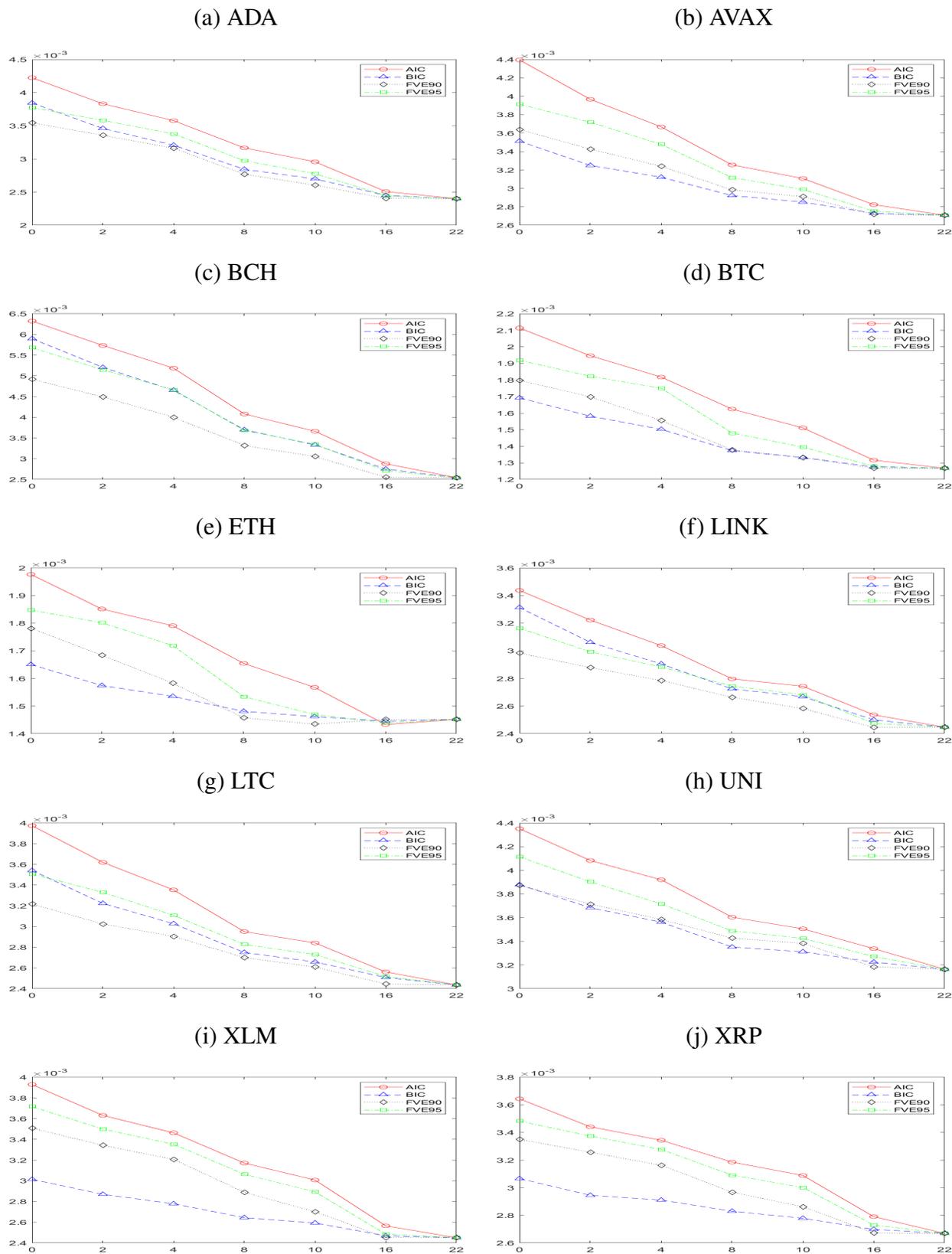

**Figure A.5** Results of the FPEs×100 based on MA with different $d$ for ten cryptocurrencies when $\tau = 0.05$.

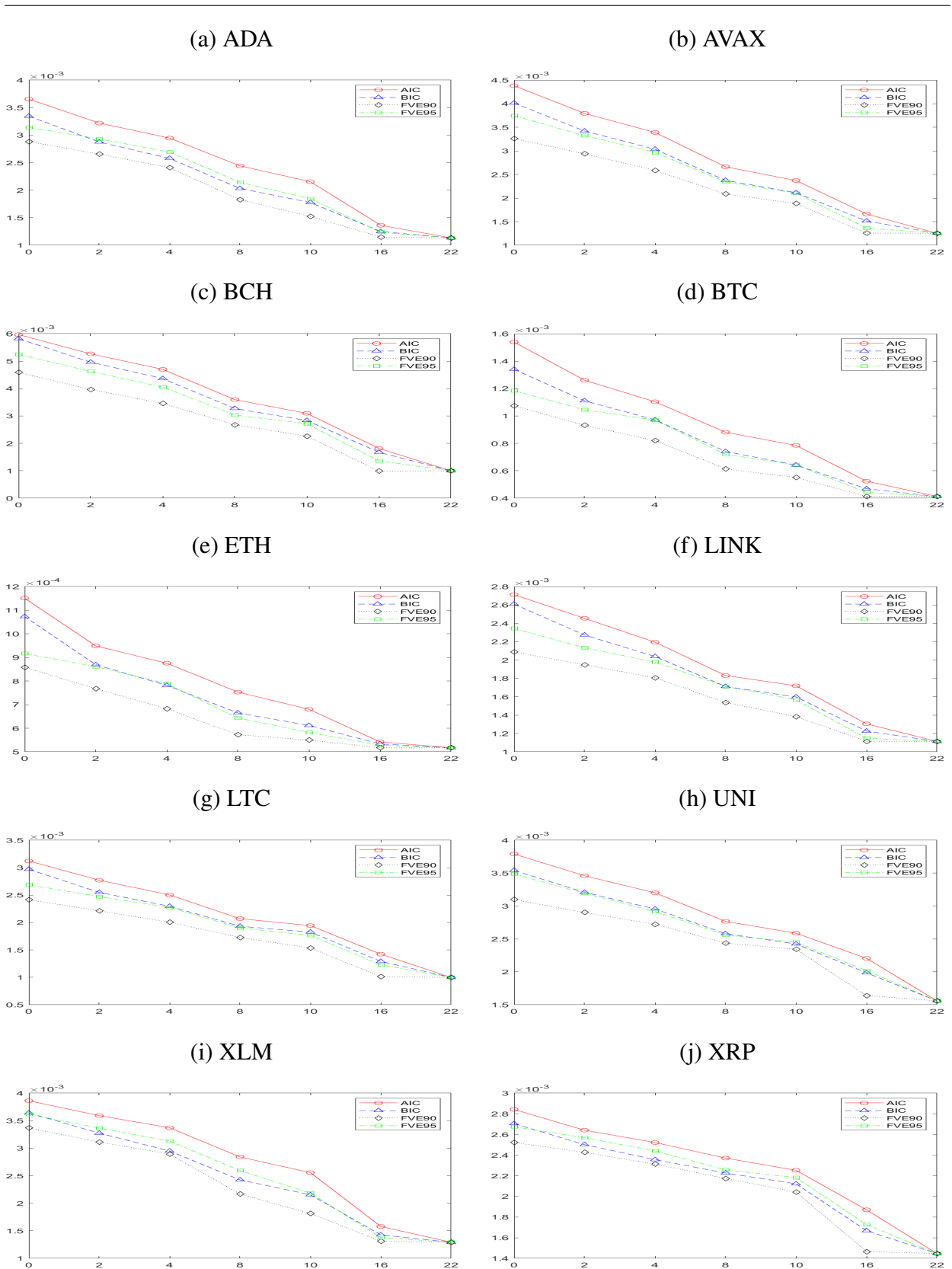

**Figure A.6** Results of the FPEs×100 based on MA with different $d$ for ten cryptocurrencies when $\tau = 0.01$.

## A2 EMPIRICAL RESULTS

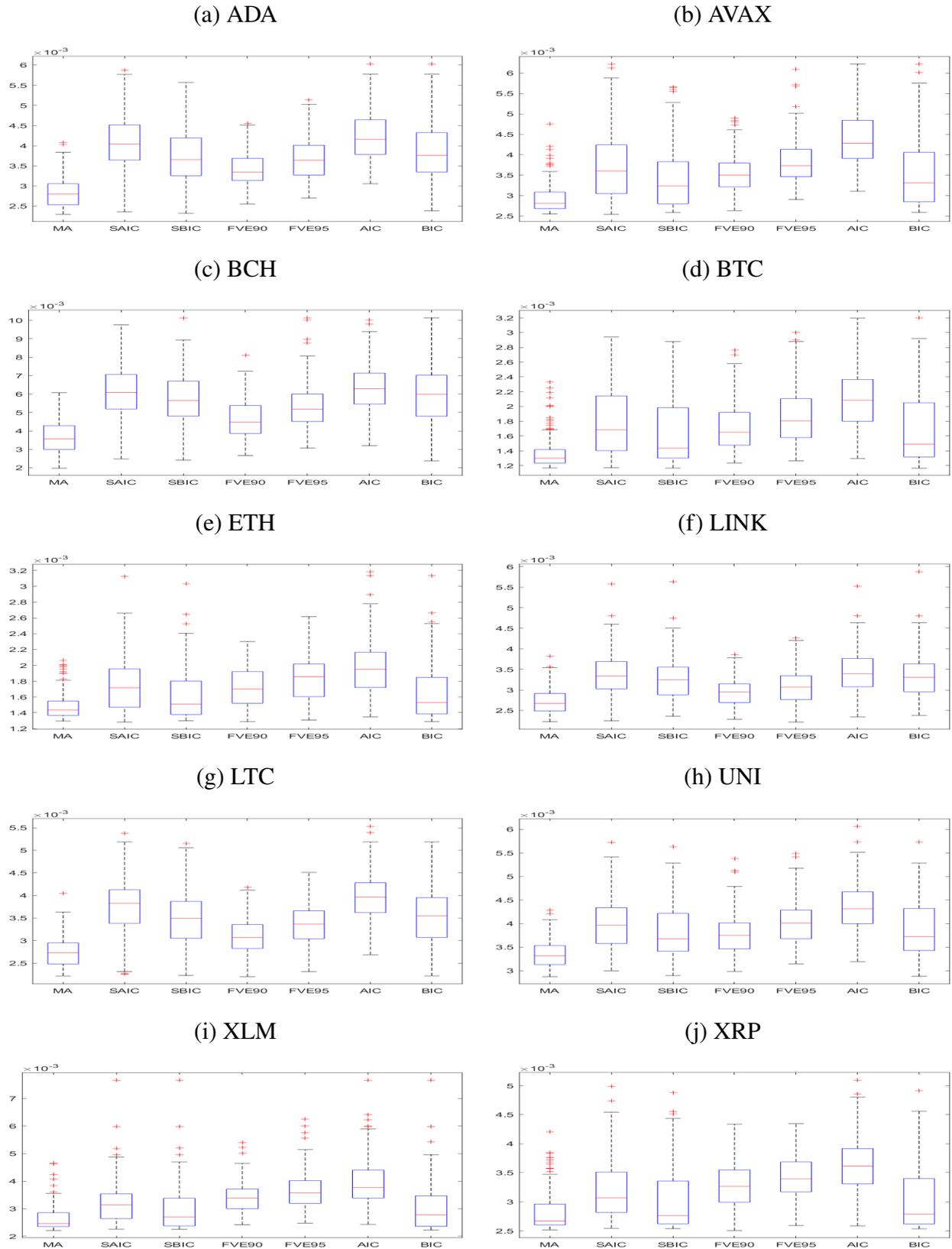

**Figure A.7** Boxplots of the FPEs×100 for ten cryptocurrencies when $\tau = 0.05$.

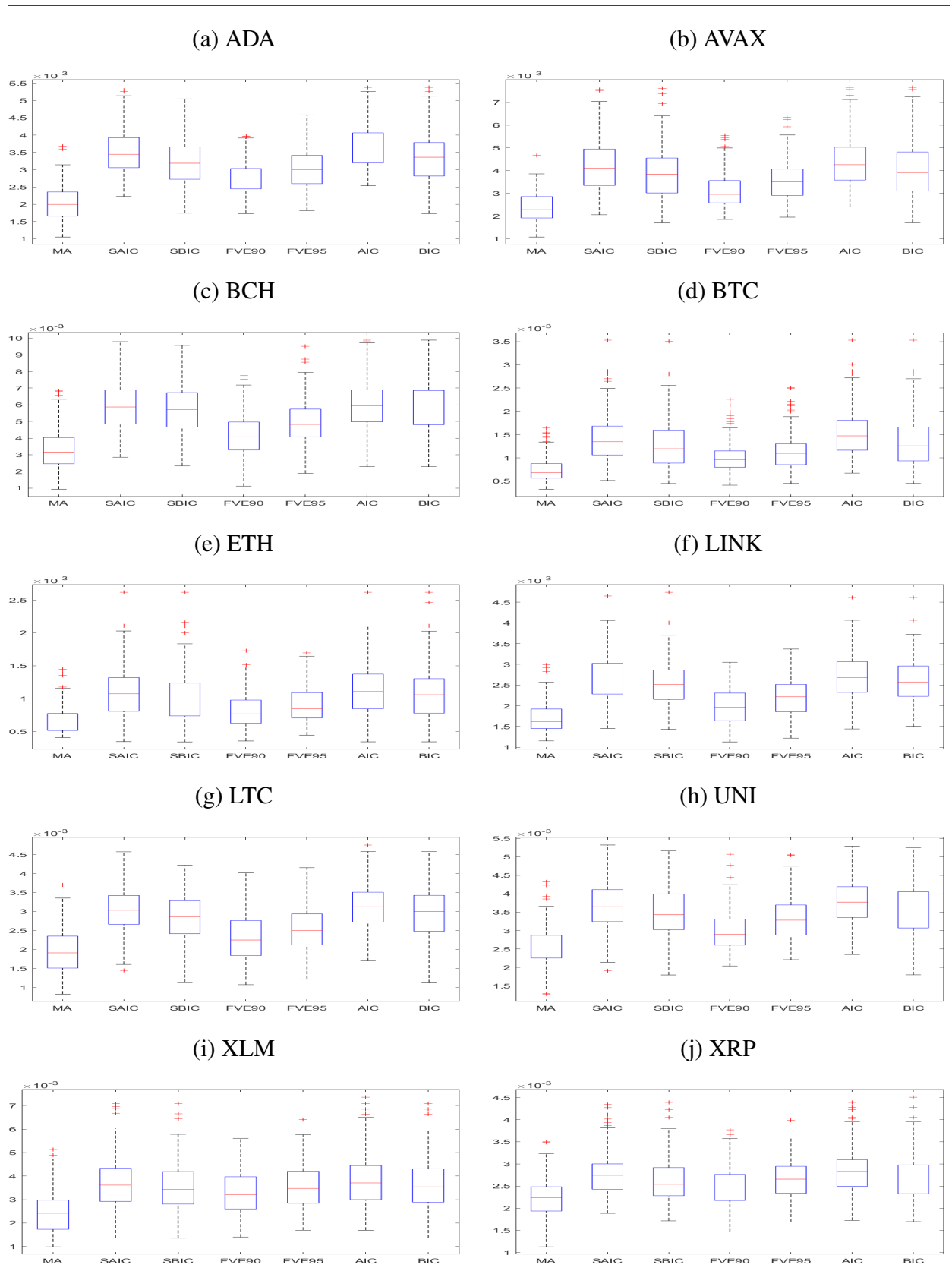

Figure A.8   Boxplots of the FPEs×100 for ten cryptocurrencies when $\tau = 0.01$.

## A2   EMPIRICAL RESULTS

**Table A.1   Results of calibration tests for the ten cryptocurrency prices**

| 5% level | MA | SAIC | SBIC | FVE90 | FVE95 | AIC | BIC |
|---|---|---|---|---|---|---|---|
| Hit | 8 | 6 | 6 | 5 | 5 | 3 | 6 |
| POF | 8 | 5 | 6 | 5 | 5 | 3 | 6 |
| CCI | 8 | 6 | 6 | 5 | 5 | 3 | 6 |
| TBF | 9 | 4 | 6 | 4 | 4 | 2 | 6 |
| 1% level | MA | SAIC | SBIC | FVE90 | FVE95 | AIC | BIC |
| Hit | 2 | 0 | 1 | 1 | 1 | 0 | 0 |
| POF | 2 | 0 | 0 | 1 | 1 | 0 | 0 |
| CCI | 2 | 0 | 1 | 1 | 1 | 0 | 1 |
| TBF | 2 | 0 | 0 | 1 | 1 | 0 | 0 |

Note: The values presented are the numbers of indices for which the test was not significantly at the 5% significance level. Larger values are better, in that they indicate the numbers of tests for which calibration was accepted.

# Supplementary Material of "Functional Data-Driven Quantile Model Averaging with Application to Cryptocurrencies"

Wenchao Xu[1], Xinyu Zhang[2], Jeng-Min Chiou[3,4] and Yuying Sun[2]

[1]*Shanghai University of International Business and Economics*

[2]*Academy of Mathematics and Systems Science, Chinese Academy of Sciences*

[3]*National Taiwan University and* [4]*Academia Sinica*

## S1. Some Details for Methodology

### S1.1. An Estimate of FPC Scores

When $X_i(t)$ are sparsely observed, $\xi_{ij}$ cannot be estimated by the integration formula. Here, the Principal components Analysis through Conditional Expectation (PACE) estimate of Yao, Müller, and Wang (2005a) is described as follows. To be more specific, let $\mathbf{U}_i = (U_{i1}, \ldots, U_{iN_i})^\top$, $\widehat{\boldsymbol{\phi}}_{ij} = \{\widehat{\phi}_j(T_{i1}), \ldots, \widehat{\phi}_j(T_{iN_i})\}^\top$, $\widehat{\boldsymbol{\mu}}_i = \{\widehat{\mu}(T_{i1}), \ldots, \widehat{\mu}(T_{iN_i})\}^\top$, and $\widehat{\boldsymbol{\Sigma}}_{\mathbf{U}_i} = \{\widehat{G}(T_{il}, T_{il'}) + \widehat{\sigma}_u^2 \delta_{ll'}\}_{l,l'=1}^{N_i}$, where $\widehat{\sigma}_u^2$ is an estimate of $\sigma_u^2$ and $\delta_{ll'} = 1$ if $l = l'$ and 0 otherwise. Then, the PACE estimator of $\xi_{ij}$ is given by $\widehat{\xi}_{ij} = \widehat{\kappa}_j \widehat{\boldsymbol{\phi}}_{ij}^\top \widehat{\boldsymbol{\Sigma}}_{\mathbf{U}_i}^{-1}(\mathbf{U}_i - \widehat{\boldsymbol{\mu}}_i)$. Note that this method can be applied to both sparsely and densely observed functional data. When $X_i(t)$ are densely observed, an alternative method is to use the "smoothing first, then estimation" procedure proposed by Ramsay and Silverman (2005), i.e., we smooth each curve first and then estimate $\xi_{ij}$ using the integration formula; see Li et al. (2010).

### S1.2. Discussions about $\text{FPE}_n(\mathbf{w})$ and $E\{\rho_\tau(\varepsilon_0)\}$

$$\begin{aligned}
&\text{FPE}_n(\mathbf{w}) - E\{\rho_\tau(\varepsilon_0)\} \\
&= E\left\{\rho_\tau\left(\varepsilon_0 + Q_\tau(X_0) - \widehat{Q}_\tau(X_0, \mathbf{w})\right) - \rho_\tau(\varepsilon_0)\Big|\mathcal{D}_n\right\} \\
&= E\left[\left\{Q_\tau(X_0) - \widehat{Q}_\tau(X_0, \mathbf{w})\right\}\psi_\tau(\varepsilon_0)\Big|\mathcal{D}_n\right] + E\left(\int_0^{\widehat{Q}_\tau(X_0,\mathbf{w})-Q_\tau(X_0)} [\mathbf{1}\{\varepsilon_0 \le s\} - \mathbf{1}\{\varepsilon_0 \le 0\}]\, ds\Big|\mathcal{D}_n\right) \\
&= E\left[\int_0^{\widehat{Q}_\tau(X_0,\mathbf{w})-Q_\tau(X_0)} \{F(s|X_0) - F(0|X_0)\}\, ds\Big|\mathcal{D}_n\right] \ge 0.
\end{aligned}$$

The last equality of Equation (S1.2) is derived by the following three facts:

$$E\left[\left\{Q_\tau(X_0) - \widehat{Q}_\tau(X_0, \mathbf{w})\right\}\psi_\tau(\varepsilon_0)\Big|\mathcal{D}_n\right]$$

## S1 SOME DETAILS FOR METHODOLOGY

$$= E\left(E\left[\left\{Q_\tau(X_0) - \widehat{Q}_\tau(X_0, \mathbf{w})\right\} \psi_\tau(\varepsilon_0) \Big| \mathcal{D}_n, \mathcal{D}_0, X_0\right] \Big| \mathcal{D}_n\right)$$

$$= E\left[\left\{Q_\tau(X_0) - \widehat{Q}_\tau(X_0, \mathbf{w})\right\} E\left\{\psi_\tau(\varepsilon_0) | \mathcal{D}_n, \mathcal{D}_0, X_0\right\} \Big| \mathcal{D}_n\right] = 0,$$

where $\mathcal{D}_0 = \{(T_{0l}, X_0(T_{0l}), \epsilon_{0l}) : l = 1, \ldots, N_0\}$ and $E\{\psi_\tau(\varepsilon_0) | \mathcal{D}_n, \mathcal{D}_0, X_0\} = E\{\psi_\tau(\varepsilon_0) | X_0\} = 0$; when $\widehat{Q}_\tau(X_0, \mathbf{w}) - Q_\tau(X_0) \geq 0$, we have

$$E\left(\int_0^{\widehat{Q}_\tau(X_0,\mathbf{w}) - Q_\tau(X_0)} [\mathbf{1}\{\varepsilon_0 \leq s\} - \mathbf{1}\{\varepsilon_0 \leq 0\}] \, ds \Big| \mathcal{D}_n\right)$$

$$= E\left(\int_0^\infty \mathbf{1}\left\{s \leq \widehat{Q}_\tau(X_0, \mathbf{w}) - Q_\tau(X_0)\right\} [\mathbf{1}\{\varepsilon_0 \leq s\} - \mathbf{1}\{\varepsilon_0 \leq 0\}] \, ds \Big| \mathcal{D}_n\right)$$

$$= \int_0^\infty E\left(\mathbf{1}\left\{s \leq \widehat{Q}_\tau(X_0, \mathbf{w}) - Q_\tau(X_0)\right\} [\mathbf{1}\{\varepsilon_0 \leq s\} - \mathbf{1}\{\varepsilon_0 \leq 0\}] \Big| \mathcal{D}_n\right) ds$$

$$= \int_0^\infty E\left\{E\left(\mathbf{1}\left\{s \leq \widehat{Q}_\tau(X_0, \mathbf{w}) - Q_\tau(X_0)\right\} [\mathbf{1}\{\varepsilon_0 \leq s\} - \mathbf{1}\{\varepsilon_0 \leq 0\}] \Big| \mathcal{D}_n, \mathcal{D}_0, X_0\right) \Big| \mathcal{D}_n\right\} ds$$

$$= \int_0^\infty E\left(\mathbf{1}\left\{s \leq \widehat{Q}_\tau(X_0, \mathbf{w}) - Q_\tau(X_0)\right\} E[\mathbf{1}\{\varepsilon_0 \leq s\} - \mathbf{1}\{\varepsilon_0 \leq 0\} | \mathcal{D}_n, \mathcal{D}_0, X_0] \Big| \mathcal{D}_n\right) ds$$

$$= \int_0^\infty E\left[\mathbf{1}\left\{s \leq \widehat{Q}_\tau(X_0, \mathbf{w}) - Q_\tau(X_0)\right\} \{F(s|X_0) - F(0|X_0)\} \Big| \mathcal{D}_n\right] ds$$

$$= E\left[\int_0^{\widehat{Q}_\tau(X_0,\mathbf{w}) - Q_\tau(X_0)} \{F(s|X_0) - F(0|X_0)\} \, ds \Big| \mathcal{D}_n\right],$$

where $E[\mathbf{1}\{\varepsilon_0 \leq s\} | \mathcal{D}_n, \mathcal{D}_0, X_0] = E[\mathbf{1}\{\varepsilon_0 \leq s\} | X_0] = F(s|X_0)$; and when $\widehat{Q}_\tau(X_0, \mathbf{w}) - Q_\tau(X_0) < 0$, the above equation (S1.2) can be derived similarly.

### S1.3. Details for the Minimization Problem (4.5)

The constrained minimization problem (4.5) can be reformulated as a linear programming problem, namely

$$\min_{\mathbf{w},\mathbf{u},\mathbf{v}} \quad \tau \mathbf{1}_n^\top \mathbf{u} + (1-\tau) \mathbf{1}_n^\top \mathbf{v}$$

$$\text{s.t.} \begin{cases} \sum_{J_n \in \mathcal{J}} w_{J_n} \widehat{Q}_{\tau,J_n}^{[-k]}(X_{(k-1)M+m}) + u_{(k-1)M+m} - v_{(k-1)M+m} = Y_{(k-1)M+m}, \\ \sum_{J_n \in \mathcal{J}} w_{J_n} = 1; w_{J_n} \geq 0 \text{ for all } J_n \in \mathcal{J}, \\ u_{(k-1)M+m} \geq 0; v_{(k-1)M+m} \geq 0 \text{ for all } k = 1, \ldots, K \text{ and } m = 1, \ldots, M, \end{cases}$$

where $\mathbf{u} = (u_1, \ldots, u_n)^\top$ and $\mathbf{v} = (v_1, \ldots, v_n)^\top$ are the positive and negative slack variables and $\mathbf{1}_n$ is the $n \times 1$ vector of ones. This linear programming can be implemented in standard software through the simplex method or the interior point method. For example, we may use the `linprog` command in MATLAB.

## S2. Conditions 1–7 for Theorems 1-2

We first introduce some conditions as follows. All limiting processes are studied with respect to $n \to \infty$.

**Condition 1** There exist a constant $a^*_{J_n}$, functions $b^*_{J_n}(t)$, $X^*_{i,J_n}(t)$, and $X^*_{0,J_n}(t)$, and a series $c_n \to 0$ such that $\widehat{a}_{J_n} - a^*_{J_n} = O_p(c_n)$, $\widehat{b}_{J_n}(t) - b^*_{J_n}(t) = O_p(c_n)$, $\widehat{\mu}(t) - \mu(t) = O_p(c_n)$, $\widehat{X}_{i,J_n}(t) - X^*_{i,J_n}(t) = O_p(c_n)$, and $E\{|\widehat{X}_{0,J_n}(t) - X^*_{0,J_n}(t)| | \mathcal{D}_n\} = O_p(c_n)$ hold uniformly for $t \in \mathcal{T}$, $J_n \in \mathcal{J}$, and $i \in \{1,\ldots,n\}$.

Let $Q^*_{\tau,J_n}(X_i) = a^*_{J_n} + \int_{\mathcal{T}} b^*_{J_n}(t)\{X^*_{i,J_n}(t) - \mu(t)\}\,dt$ and $Q^*_\tau(X_i, \mathbf{w}) = \sum_{J_n \in \mathcal{J}} w_{J_n} Q^*_{\tau,J_n}(X_i)$, where $i = 0, 1, \ldots, n$. Define $\text{EFPE}^*_n(\mathbf{w}) = E\{\rho_\tau(Y_0 - Q^*_\tau(X_0, \mathbf{w}))\} - E\{\rho_\tau(\varepsilon_0)\} = E\left[\int_0^{Q^*_\tau(X_0,\mathbf{w}) - Q_\tau(X_0)}\{F(s|X_0) - F(0|X_0)\}\,ds\right] \geq 0$, where the second equality is due to Knight's identity (4.4). Let $\eta_n = n \inf_{\mathbf{w} \in \mathcal{W}} \text{EFPE}^*_n(\mathbf{w})$.

**Condition 2** There exists a constant $\varrho > 0$ such that $|Q^*_{\tau,J_n}(X_i) - Q_\tau(X_i)| \leq \varrho$ a.s. uniformly for $i \in \{0, 1, \ldots, n\}$ and $J_n \in \mathcal{J}$.

**Condition 3** There exists a series $g_n \to 0$ such that $\widehat{Q}_\tau(X_{(k-1)M+m}, \mathbf{w}) - \widehat{Q}^{[-k]}_\tau(X_{(k-1)M+m}, \mathbf{w}) = O_p(g_n)$ holds uniformly for $\mathbf{w} \in \mathcal{W}$, $k \in \{1,\ldots,K\}$, and $m \in \{1,\ldots,M\}$.

**Condition 4** $\eta_n^{-1} n^{1/2} |\mathcal{J}|^{1/2} = o(1)$, $\eta_n^{-1} n c_n = o(1)$, and $\eta_n^{-1} n g_n = o(1)$.

Condition 1 describes the convergence rates of the estimators $\widehat{a}_{J_n}$, $\widehat{b}_{J_n}(t)$, $\widehat{\mu}(t)$, and $\widehat{X}_{i,J_n}(t)$ under each model. They need not have the same convergence rate. For example, when $\widehat{a}_{J_n} - a^*_{J_n} = O_p(n^{-\alpha_1})$, $\widehat{b}_{J_n}(t) - b^*_{J_n}(t) = O_p(n^{-\alpha_2})$, $\widehat{\mu}(t) - \mu(t) = O_p(n^{-\alpha_3})$, $\widehat{X}_{i,J_n}(t) - X^*_{i,J_n}(t) = O_p(n^{-\alpha_4})$, and $E\{|\widehat{X}_{0,J_n}(t) - X^*_{0,J_n}(t)| | \mathcal{D}_n\} = O_p(n^{-\alpha_5})$, $c_n$ can be $n^{-\min\{\alpha_1,\ldots,\alpha_5\}}$. This condition is very mild and is verified in the functional data literature; see, e.g., Yao et al. (2005a,b), Kato (2012), and Li et al. (2022). Note that when model $J_n$ is a true model, $a^*_{J_n}$, $b^*_{J_n}(t)$, $X^*_{i,J_n}(t)$, and $X^*_{0,J_n}(t)$ are naturally the true parameter values. Condition 2 excludes some pathological cases in which $Q^*_{\tau,J_n}(X_i) - Q_\tau(X_i)$ explodes. Condition 3 is similar to the condition 6 of Zhang et al. (2018), which requires the difference between the regular prediction $\widehat{Q}_\tau(X_{(k-1)M+m}, \mathbf{w})$ and the leave-$M$-out prediction $\widehat{Q}^{[-k]}_\tau(X_{(k-1)M+m}, \mathbf{w})$ to decrease with the rate $g_n$ as the sample size increases. Since $\widehat{Q}_\tau(X_{(k-1)M+m}, \mathbf{w}) - \widehat{Q}^{[-k]}_\tau(X_{(k-1)M+m}, \mathbf{w})$

$$= \sum_{J_n \in \mathcal{J}} w_{J_n} \left\{ \widehat{Q}_{\tau, J_n}(X_{(k-1)M+m}) - \widehat{Q}_{\tau, J_n}^{[-k]}(X_{(k-1)M+m}) \right\}$$

$$= \sum_{J_n \in \mathcal{J}} w_{J_n} \left( \widehat{a}_{J_n} - \widehat{a}_{J_n}^{[-k]} \right) + \sum_{J_n \in \mathcal{J}} w_{J_n} \int_{\mathcal{T}} \left\{ \widehat{b}_{J_n}(t) - \widehat{b}_{J_n}^{[-k]}(t) \right\} \left\{ \widehat{X}_{(k-1)M+m, J_n}(t) - \widehat{\mu}(t) \right\} dt,$$

the explicit expression of $g_n$ can be derived by analyzing the orders of $\max_{J_n, k} |\widehat{a}_{J_n} - \widehat{a}_{J_n}^{[-k]}|$ and $\sup_{t \in \mathcal{T}} \max_{J_n, k} |\widehat{b}_{J_n}(t) - \widehat{b}_{J_n}^{[-k]}(t)|$. Condition 4 requires that $\eta_n$ grows at a rate no slower than $\max\{n^{1/2}|\mathcal{J}|^{1/2}, nc_n, ng_n\}$. Note that when there is at least one true model included in the set of candidate models, $\eta_n \equiv 0$. Therefore, Condition 4 implies that all candidate models are misspecified. This condition is similar to Condition (21) of Zhang, Wan, and Zou (2013), Condition (7) of Ando and Li (2014), and Condition 3 of Zhang, Chiou, and Ma (2018). The first part of Condition 4 also implies that the number of candidate models $|\mathcal{J}|$ is allowed to grow to infinity as the sample size increases. Condition 4 excludes the lucky situation where one of these candidate models happens to be true. After that, we demonstrate that when there is at least one true model included in the set of candidate models, our averaged regression estimation is consistent.

Let $f(\cdot|X_i)$ denote the conditional density function of $\varepsilon_i$ given $X_i$, and let $\widehat{X}_i(t) = \widehat{\mu}(t) + \sum_{j=1}^{\infty} \widehat{\xi}_{ij} \widehat{\phi}_j(t)$, $i = 1, \ldots, n$. We further impose the following additional conditions for Theorem 2.

**Condition 5** There exist constants $0 < \underline{C} \leq \bar{C} < \infty$ such that $\underline{C} \leq f(s|X_i) \leq \bar{C}$ holds uniformly for $|s| \leq \varrho$ and $i = 1, \ldots, n$, where $\varrho$ is defined in Condition 2.

**Condition 6** $\widehat{X}_i(t)$ and $\widehat{\mu}(t)$ are $O_p(1)$ uniformly for $i \in \{1, \ldots, n\}$ and $t \in \mathcal{T}$.

**Condition 7** There exists $\lambda_{\min} > 0$ such that for uniformly for $\mathbf{w} \in \mathcal{W}$ and for almost all $i \in \{1, \ldots, n\}$, $\left[ \widehat{a}_{\mathbf{w}} - a + \int_{\mathcal{T}} \left\{ \widehat{b}_{\mathbf{w}}(t) - b(t) \right\} \left\{ \widehat{X}_i(t) - \widehat{\mu}(t) \right\} dt \right]^2 \geq \lambda_{\min} \left[ (\widehat{a}_{\mathbf{w}} - a)^2 + \int_{\mathcal{T}} \left\{ \widehat{b}_{\mathbf{w}}(t) - b(t) \right\}^2 dt \right]$.

Condition 5 is mild, and it allows conditional heteroskedasticity in the FLQR model. Condition 6 is also mild and is the same as Condition 8 of Zhang et al. (2018), which excludes some pathological cases where $\widehat{X}_i(t)$ and $\widehat{\mu}(t)$ explode. Condition 7 is similar to Condition 7 of Zhang et al. (2018), which states that most $\widehat{X}_i(t)$ do not degenerate in the sense that their inner products with $\widehat{b}_{\mathbf{w}}(t) - b(t)$ do not approach zero.

## S3. Proofs of Theorems 1 and 2

### S3.1. Technical Lemmas

Define the metric $\|\cdot\|_1$ on $\mathcal{W}$ as

$$\|\mathbf{w}_1 - \mathbf{w}_2\|_1 = \sum_{J_n \in \mathcal{J}} |w_{1,J_n} - w_{2,J_n}|$$

for any $\mathbf{w}_1 = \{w_{1,J_n} : J_n \in \mathcal{J}\}$ and $\mathbf{w}_2 = \{w_{2,J_n} : J_n \in \mathcal{J}\} \in \mathcal{W}$. Denote $\mathcal{N}(\epsilon, \mathcal{W}, \|\cdot\|_1)$ as the $\epsilon$-covering number of $\mathcal{W}$. By the proposition C.1 of Ghosal and van der Vaart (2017), we have for $0 < \epsilon \le 1$,

$$\mathcal{N}(\epsilon, \mathcal{W}, \|\cdot\|_1) \le \left(1 + \frac{4}{\epsilon}\right)^{|\mathcal{J}|-1}. \tag{S3.-3}$$

Let $\xi_i$, $i = 1, \ldots, n$ be an i.i.d. sequence of random variables or vectors, $P_n$ denote the empirical distribution based on them, and $\mathcal{H}$ be a class of functions defined on them. The $L_2(P_n)$-norm of $\mathcal{H}$ is defined by $\|h\|_{P_n} = \sqrt{n^{-1} \sum_{i=1}^n h(\xi_i)^2}$ for $h \in \mathcal{H}$. We first present a maximal inequality in the following lemma.

**Lemma 1.** *Assume that $\sup_{h \in \mathcal{H}} \|h\|_{P_n} \le \varrho_1$ for some constant $\varrho_1 > 0$. Let $W_i$, $i = 1, \ldots, n$ be mean zero independent random variables such that $|W_i| \le B$ a.s. for some constant $B > 0$. Then, there exists a constant $C^* > 0$ depending only on $B$ such that for all $n \ge \delta^{*-2} C^{*2} \varrho_1^2$ satisfying*

$$\frac{C^*}{\sqrt{2n}} \int_{\delta^*/8B}^{\varrho_1} \sqrt{\log \mathcal{N}(u, \mathcal{H}, \|\cdot\|_{P_n})} \, du \le \delta^*,$$

*we have*

$$\mathbb{P}\left(\sup_{h \in \mathcal{H}} \left|\frac{1}{n} \sum_{i=1}^n W_i h(\xi_i)\right| \ge \delta^*\right) \le C^* \exp\left(-\frac{n \delta^{*2}}{C^{*2} \varrho_1^2}\right).$$

Lemma 1 follows easily from the lemma 3.2 of van de Geer (2000). Applying (S3.-3) and Lemma 1, we obtain the following main lemma.

**Lemma 2.** *Let $h(\cdot, \cdot; \mathbf{w}) \colon L_2(\mathcal{T}) \times \mathbb{R} \to \mathbb{R}$ be a measurable function for $\mathbf{w} \in \mathcal{W}$. Assume $|h(\cdot, \cdot; \mathbf{w})| \le \varrho_1$ and $|h(\cdot, \cdot; \mathbf{w}_1) - h(\cdot, \cdot; \mathbf{w}_2)| \le \varrho_1 \|\mathbf{w}_1 - \mathbf{w}_2\|_1$ for some constant $\varrho_1 > 0$ and any $\mathbf{w}_1, \mathbf{w}_2 \in \mathcal{W}$. Then,*

$$\sup_{\mathbf{w} \in \mathcal{W}} \left|\frac{1}{n} \sum_{i=1}^n h(X_i, \varepsilon_i; \mathbf{w}) - E\{h(X_i, \varepsilon_i; \mathbf{w})\}\right| = O_p(n^{-1/2} |\mathcal{J}|^{1/2}).$$

## S3 PROOFS OF THEOREMS 1 AND 2

*Proof of Lemma 2.* Let $\mathcal{H} = \{h(\cdot, \cdot; \mathbf{w}) : \mathbf{w} \in \mathcal{W}\}$ and $P_n$ denote the empirical distribution based on $\{(X_i, \varepsilon_i) : i = 1, \ldots, n\}$. Then,

$$\sup_{\mathbf{w} \in \mathcal{W}} \left| \frac{1}{n} \sum_{i=1}^n h(X_i, \varepsilon_i; \mathbf{w}) - E\{h(X_i, \varepsilon_i; \mathbf{w})\} \right| = \sup_{h \in \mathcal{H}} \left| \int h \, d(P_n - P) \right|. \tag{S3.-2}$$

Let $W_1, \ldots, W_n$ be i.i.d. Rademacher random variables independent of $\{(X_i, \varepsilon_i) : i = 1, \ldots, n\}$, i.e., $\mathbb{P}(W_i = 1) = \mathbb{P}(W_i = -1) = 1/2$. By the corollary 3.4 of van de Geer (2000), for $n \geq 8\varrho_1^2/\delta^2$,

$$\mathbb{P}\left( \sup_{h \in \mathcal{H}} \left| \int h \, d(P_n - P) \right| \geq \delta \right) \leq 4\mathbb{P}\left( \sup_{h \in \mathcal{H}} \left| \frac{1}{n} \sum_{i=1}^n W_i h(X_i, \epsilon_i) \right| \geq \delta/4 \right). \tag{S3.-1}$$

The next task is to bound the right term of (S3.-1). Using the condition $|h(\cdot, \cdot; \mathbf{w}_1) - h(\cdot, \cdot; \mathbf{w}_2)| \leq \varrho_1 \|\mathbf{w}_1 - \mathbf{w}_2\|_1$ for any $\mathbf{w}_1, \mathbf{w}_2 \in \mathcal{W}$, we have

$$\|h_1 - h_2\|_{P_n} = \left[ \frac{1}{n} \sum_{i=1}^n \{h(X_i, \epsilon_i; \mathbf{w}_1) - h(X_i, \epsilon_i; \mathbf{w}_2)\}^2 \right]^{1/2} \leq \varrho_1 \|\mathbf{w}_1 - \mathbf{w}_2\|_1$$

for any $h_1 = h(\cdot, \cdot, \mathbf{w}_1)$ and $h_2 = h(\cdot, \cdot, \mathbf{w}_2) \in \mathcal{H}$. Therefore, from (S3.-3), $\mathcal{N}(u, \mathcal{H}, \|\cdot\|_{P_n}) \leq (1 + 4\varrho_1/u)^{|\mathcal{J}|-1}$. Thus, we have

$$\int_0^{\varrho_1} \sqrt{\log \mathcal{N}(u, \mathcal{H}, \|\cdot\|_{P_n})} \, du \leq |\mathcal{J}|^{1/2} \int_0^{\varrho_1} \sqrt{\log(1 + 4\varrho_1/u)} \, du$$

$$\leq |\mathcal{J}|^{1/2} \int_0^{\varrho_1} \log(1 + 4\varrho_1/u) \, du = |\mathcal{J}|^{1/2} \varrho_1 \int_0^1 \log(1 + 4/u) \, du$$

$$= |\mathcal{J}|^{1/2} \varrho_1 \int_0^1 [\{\log(u+4) + 1\} - (\log u + 1)] \, du$$

$$= |\mathcal{J}|^{1/2} \varrho_1 (5 \log 5 - 4 \log 4)$$

$$\leq 3 |\mathcal{J}|^{1/2} \varrho_1,$$

where the second inequality follows from $1 + 4\varrho_1/u \geq e$ for all $0 < u \leq \varrho_1$. By taking $\delta^* = Cn^{-1/2}|\mathcal{J}|^{1/2}$ in Lemma 1 for $C \geq C^* \varrho_1 \max\{3/\sqrt{2}, |\mathcal{J}|^{-1/2}\}$, we have

$$\mathbb{P}\left( \sup_{h \in \mathcal{H}} \left| \frac{1}{n} \sum_{i=1}^n W_i h(X_i, \epsilon_i) \right| \geq Cn^{-1/2}|\mathcal{J}|^{1/2} \right) \leq C^* \exp\left( -\frac{C^2 |\mathcal{J}|}{C^{*2} \varrho_1^2} \right).$$

Therefore, with $\delta = 4Cn^{-1/2}|\mathcal{J}|^{1/2}$ in (S3.-1) for $C \geq \varrho_1 \max\{3C^*/\sqrt{2}, C^*|\mathcal{J}|^{-1/2}, |\mathcal{J}|^{-1/2}/\sqrt{2}\}$, we obtain

$$\mathbb{P}\left(\sup_{h \in \mathcal{H}} \left|\int h \, d(P_n - P)\right| \geq Cn^{-1/2}|\mathcal{J}|^{1/2}\right) \leq 4C^* \exp\left(-\frac{C^2|\mathcal{J}|}{C^{*2}\varrho_1^2}\right),$$

which yields

$$\sup_{h \in \mathcal{H}} \left|\int h \, d(P_n - P)\right| = O_p(n^{-1/2}|\mathcal{J}|^{1/2}). \tag{S3.0}$$

Combining (S3.-2) and (S3.0), the proof of Lemma 2 is completed.

### S3.2. Proof of Theorem 1

Let $\widetilde{\text{CV}}_K(\mathbf{w}) = \text{CV}_K(\mathbf{w}) - n^{-1}\sum_{i=1}^{n} \rho_\tau(\varepsilon_i)$, where the second term is unrelated to $\mathbf{w}$. Therefore, $\widehat{\mathbf{w}} = \arg\min_{\mathbf{w} \in \mathcal{W}} \widetilde{\text{CV}}_K(\mathbf{w})$. To prove Theorem 1, it is sufficient to show that

$$\sup_{\mathbf{w} \in \mathcal{W}} \frac{|\text{EFPE}_n(\mathbf{w}) - \text{EFPE}_n^*(\mathbf{w})|}{\text{EFPE}_n^*(\mathbf{w})} = O_p(\eta_n^{-1} nc_n) \tag{S3.1}$$

and

$$\sup_{\mathbf{w} \in \mathcal{W}} \frac{\left|\widetilde{\text{CV}}_K(\mathbf{w}) - \text{EFPE}_n^*(\mathbf{w})\right|}{\text{EFPE}_n^*(\mathbf{w})} = O_p\left\{\eta_n^{-1}\left(nc_n + ng_n + n^{1/2}|\mathcal{J}|^{1/2}\right)\right\}. \tag{S3.2}$$

The reason is as follows. First, Condition 2 implies that $\text{EFPE}_n^*(\mathbf{w}) \leq \varrho$ for any $\mathbf{w} \in \mathcal{W}$, which along with Condition 4 and (S3.1) leads to $\inf_{\mathbf{w} \in \mathcal{W}} \text{EFPE}_n(\mathbf{w}) \leq \varrho\{1 + o_p(1)\}$, where $\varrho$ is defined in Condition 2. From the definition of the infimum, there exist a non-negative series $\vartheta_n \to 0$ and a vector $\mathbf{w}_n^* \in \mathcal{W}$ such that $\inf_{\mathbf{w} \in \mathcal{W}} \text{EFPE}_n(\mathbf{w}) + \vartheta_n = \text{EFPE}_n(\mathbf{w}_n^*)$. Observe that $\widetilde{\text{CV}}_K(\widehat{\mathbf{w}}) \leq \widetilde{\text{CV}}_K(\mathbf{w}_n^*)$. Therefore,

$$\frac{\text{EFPE}_n(\widehat{\mathbf{w}})}{\text{EFPE}_n(\mathbf{w}_n^*)}\left(\frac{\widetilde{\text{CV}}_K(\widehat{\mathbf{w}}) - \text{EFPE}_n(\widehat{\mathbf{w}})}{\text{EFPE}_n(\widehat{\mathbf{w}})} + 1\right) \leq \frac{\widetilde{\text{CV}}_K(\mathbf{w}_n^*) - \text{EFPE}_n(\mathbf{w}_n^*)}{\text{EFPE}_n(\mathbf{w}_n^*)} + 1. \tag{S3.3}$$

By (S3.1), (S3.2), and Condition 4, we have

$$\sup_{\mathbf{w} \in \mathcal{W}} \frac{\left|\widetilde{\text{CV}}_K(\mathbf{w}) - \text{EFPE}_n(\mathbf{w})\right|}{\text{EFPE}_n(\mathbf{w})}$$

$$\leq \left(\sup_{\mathbf{w} \in \mathcal{W}} \frac{\left|\widetilde{\text{CV}}_K(\mathbf{w}) - \text{EFPE}_n^*(\mathbf{w})\right|}{\text{EFPE}_n^*(\mathbf{w})} + \sup_{\mathbf{w} \in \mathcal{W}} \frac{|\text{EFPE}_n(\mathbf{w}) - \text{EFPE}_n^*(\mathbf{w})|}{\text{EFPE}_n^*(\mathbf{w})}\right) \sup_{\mathbf{w} \in \mathcal{W}} \frac{\text{EFPE}_n^*(\mathbf{w})}{\text{EFPE}_n(\mathbf{w})}$$

$$= O_p\{\eta_n^{-1}\left(nc_n + ng_n + n^{1/2}|\mathcal{J}|^{1/2}\right)\},$$

## S3 PROOFS OF THEOREMS 1 AND 2

where the last step is due to $\sup_{\mathbf{w}\in\mathcal{W}} \frac{\text{EFPE}_n^*(\mathbf{w})}{\text{EFPE}_n(\mathbf{w})} = O_p(1)$. Therefore, from (S3.3), we have

$$\frac{\text{EFPE}_n(\widehat{\mathbf{w}})}{\text{EFPE}_n(\mathbf{w}_n^*)} - 1 \leq \left( \frac{\widetilde{\text{CV}}_K(\mathbf{w}_n^*) - \text{EFPE}_n(\mathbf{w}_n^*)}{\text{EFPE}_n(\mathbf{w}_n^*)} - \frac{\widetilde{\text{CV}}_K(\widehat{\mathbf{w}}) - \text{EFPE}_n(\widehat{\mathbf{w}})}{\text{EFPE}_n(\widehat{\mathbf{w}})} \right)$$
$$\times \left( \frac{\widetilde{\text{CV}}_K(\widehat{\mathbf{w}}) - \text{EFPE}_n(\widehat{\mathbf{w}})}{\text{EFPE}_n(\widehat{\mathbf{w}})} + 1 \right)^{-1}$$
$$= O_p\{\eta_n^{-1}\left(nc_n + ng_n + n^{1/2}|\mathcal{J}|^{1/2}\right)\},$$

which implies that

$$0 \leq \frac{\text{EFPE}_n(\widehat{\mathbf{w}})}{\inf_{\mathbf{w}\in\mathcal{W}} \text{EFPE}_n(\mathbf{w})} - 1 = \frac{\text{EFPE}_n(\widehat{\mathbf{w}})}{\text{EFPE}_n(\mathbf{w}_n^*)} \left( 1 + \frac{\vartheta_n}{\inf_{\mathbf{w}\in\mathcal{W}} \text{EFPE}_n(\mathbf{w})} \right) - 1$$
$$\leq \left[ 1 + O_p\left\{\eta_n^{-1}\left(nc_n + ng_n + n^{1/2}|\mathcal{J}|^{1/2}\right)\right\} \right] \left( 1 + \frac{\vartheta_n}{\inf_{\mathbf{w}\in\mathcal{W}} \text{EFPE}_n(\mathbf{w})} \right) - 1.$$

By taking $\vartheta_n/\inf_{\mathbf{w}\in\mathcal{W}} \text{EFPE}_n(\mathbf{w}) = O_p\{\eta_n^{-1}\left(nc_n + ng_n + n^{1/2}|\mathcal{J}|^{1/2}\right)\}$, Theorem 1 holds. The remainder of this proof is to show (S3.1) and (S3.2). Using Condition 1, we obtain

$$E\left\{\left|\widehat{Q}_{\tau,J_n}(X_0) - Q^*_{\tau,J_n}(X_0)\right|\Big|\mathcal{D}_n\right\}$$
$$\leq |\widehat{a}_{J_n} - a^*_{J_n}| + \left|\int_{\mathcal{T}} \left\{\widehat{b}_{J_n}(t) - b^*_{J_n}(t)\right\} \left\{X^*_{0,J_n}(t) - \mu(t)\right\} dt\right|$$
$$+ \int_{\mathcal{T}} \left|\widehat{b}_{J_n}(t)\right| \left[E\left\{\left|\widehat{X}_{0,J_n}(t) - X^*_{0,J_n}(t)\right|\Big|\mathcal{D}_n\right\} - |\widehat{\mu}(t) - \mu(t)|\right] dt$$
$$= O_p(c_n)$$

uniformly for $J_n \in \mathcal{J}$, which, along with the fact that $|F(s|X_0) - F(0|X_0)| \leq 1$, implies

$$\sup_{\mathbf{w} \in \mathcal{W}} |\text{EFPE}_n(\mathbf{w}) - \text{EFPE}_n^*(\mathbf{w})|/\text{EFPE}_n^*(\mathbf{w})$$

$$\leq \eta_n^{-1} n \sup_{\mathbf{w} \in \mathcal{W}} |\text{EFPE}_n(\mathbf{w}) - \text{EFPE}_n^*(\mathbf{w})|$$

$$= \eta_n^{-1} n \sup_{\mathbf{w} \in \mathcal{W}} \left| E\left[ \int_{Q_\tau^*(X_0,\mathbf{w}) - Q_\tau(X_0)}^{\widehat{Q}_\tau(X_0,\mathbf{w}) - Q_\tau(X_0)} \{F(s|X_0) - F(0|X_0)\} \, ds \Big| \mathcal{D}_n \right] \right|$$

$$\leq \eta_n^{-1} n \sup_{\mathbf{w} \in \mathcal{W}} E\left\{ \left\| \widehat{Q}_\tau(X_0, \mathbf{w}) - Q_\tau^*(X_0, \mathbf{w}) \right\| \Big| \mathcal{D}_n \right\}$$

$$\leq \eta_n^{-1} n \sup_{\mathbf{w} \in \mathcal{W}} \sum_{J_n \in \mathcal{J}} w_{J_n} E\left\{ \left\| \widehat{Q}_{\tau, J_n}(X_0) - Q_{\tau, J_n}^*(X_0) \right\| \Big| \mathcal{D}_n \right\}$$

$$\leq \eta_n^{-1} n \max_{J_n \in \mathcal{J}} E\left\{ \left\| \widehat{Q}_{\tau, J_n}(X_0) - Q_{\tau, J_n}^*(X_0) \right\| \Big| \mathcal{D}_n \right\}$$

$$= O_p(\eta_n^{-1} n c_n).$$

Finally, (S3.1) is obtained. By Knight's identity (4.4), we have

$$\widetilde{\text{CV}}_K(\mathbf{w}) - \text{EFPE}_n^*(\mathbf{w})$$

$$= \frac{1}{n} \sum_{k=1}^{K} \sum_{m=1}^{M} \left\{ \rho_\tau\left( Y_{(k-1)M+m} - \widehat{Q}_\tau^{[-k]}(X_{(k-1)M+m}, \mathbf{w}) \right) - \rho_\tau(\varepsilon_{(k-1)M+m}) \right\} - \text{EFPE}_n^*(\mathbf{w})$$

$$= \frac{1}{n} \sum_{k=1}^{K} \sum_{m=1}^{M} \left\{ Q_\tau(X_{(k-1)M+m}) - \widehat{Q}_\tau^{[-k]}(X_{(k-1)M+m}, \mathbf{w}) \right\} \psi_\tau(\varepsilon_{(k-1)M+m})$$

$$+ \frac{1}{n} \sum_{k=1}^{K} \sum_{m=1}^{M} \int_0^{\widehat{Q}_\tau^{[-k]}(X_{(k-1)M+m},\mathbf{w}) - Q_\tau(X_{(k-1)M+m})} \left[ \mathbf{1}\{\varepsilon_{(k-1)M+m} \leq s\} - \mathbf{1}\{\varepsilon_{(k-1)M+m} \leq 0\} \right] ds$$

$$- \text{EFPE}_n^*(\mathbf{w})$$

$$\equiv \text{CV}_{1n}(\mathbf{w}) + \text{CV}_{2n}(\mathbf{w}) + \text{CV}_{3n}(\mathbf{w}) + \text{CV}_{4n}(\mathbf{w}),$$

where

$$\text{CV}_{1n}(\mathbf{w}) = \frac{1}{n} \sum_{i=1}^{n} \{Q_\tau(X_i) - Q_\tau^*(X_i, \mathbf{w})\} \psi_\tau(\varepsilon_i),$$

$$\text{CV}_{2n}(\mathbf{w}) = \frac{1}{n} \sum_{i=1}^{n} \int_0^{Q_\tau^*(X_i,\mathbf{w}) - Q_\tau(X_i)} \left[ \mathbf{1}\{\varepsilon_i \leq s\} - \mathbf{1}\{\varepsilon_i \leq 0\} \right] ds - \text{EFPE}_n^*(\mathbf{w}),$$

$$\text{CV}_{3n}(\mathbf{w}) = \frac{1}{n} \sum_{i=1}^{n} \left\{ Q_\tau^*(X_i, \mathbf{w}) - \widehat{Q}_\tau(X_i, \mathbf{w}) \right\} \psi_\tau(\varepsilon_i)$$

$$+ \frac{1}{n} \sum_{i=1}^{n} \int_{Q_\tau^*(X_i,\mathbf{w}) - Q_\tau(X_i)}^{\widehat{Q}_\tau(X_i,\mathbf{w}) - Q_\tau(X_i)} \left[ \mathbf{1}\{\varepsilon_i \leq s\} - \mathbf{1}\{\varepsilon_i \leq 0\} \right] ds,$$

# S3 PROOFS OF THEOREMS 1 AND 2

and

$$\text{CV}_{4n}(\mathbf{w}) = \frac{1}{n} \sum_{k=1}^{K} \sum_{m=1}^{M} \left\{ \widehat{Q}_\tau(X_{(k-1)M+m}, \mathbf{w}) - \widehat{Q}_\tau^{[-k]}(X_{(k-1)M+m}, \mathbf{w}) \right\} \psi_\tau(\varepsilon_{(k-1)M+m})$$

$$+ \frac{1}{n} \sum_{k=1}^{K} \sum_{m=1}^{M} \int_{\widehat{Q}_\tau(X_{(k-1)M+m},\mathbf{w}) - Q_\tau(X_{(k-1)M+m})}^{\widehat{Q}_\tau^{[-k]}(X_{(k-1)M+m},\mathbf{w}) - Q_\tau(X_{(k-1)M+m})} \left[ \mathbf{1}\{\varepsilon_{(k-1)M+m} \leq s\} - \mathbf{1}\{\varepsilon_{(k-1)M+m} \leq 0\} \right] ds.$$

Hence, to prove (S3.2), it is sufficient to verify that

$$\sup_{\mathbf{w} \in \mathcal{W}} |\text{CV}_{1n}(\mathbf{w})|/\text{EFPE}_n^*(\mathbf{w}) = O_p(\eta_n^{-1} n^{1/2} |\mathcal{J}|^{1/2}), \tag{S3.3}$$

$$\sup_{\mathbf{w} \in \mathcal{W}} |\text{CV}_{2n}(\mathbf{w})|/\text{EFPE}_n^*(\mathbf{w}) = O_p(\eta_n^{-1} n^{1/2} |\mathcal{J}|^{1/2}), \tag{S3.4}$$

$$\sup_{\mathbf{w} \in \mathcal{W}} |\text{CV}_{3n}(\mathbf{w})|/\text{EFPE}_n^*(\mathbf{w}) = O_p(\eta_n^{-1} n c_n), \tag{S3.5}$$

and

$$\sup_{\mathbf{w} \in \mathcal{W}} |\text{CV}_{4n}(\mathbf{w})|/\text{EFPE}_n^*(\mathbf{w}) = O_p(\eta_n^{-1} n g_n). \tag{S3.6}$$

*Proof of* (S3.3). Take

$$h(X_i, \epsilon_i; \mathbf{w}) = \{Q_\tau(X_i) - Q_\tau^*(X_i, \mathbf{w})\} \psi_\tau(\varepsilon_i),$$

in Lemma 2. By the fact $|\psi_\tau(\varepsilon_i)| \leq 1$ and Condition 2, it is easy to see that $|h(X_i, \epsilon_i; \mathbf{w})| \leq \varrho$ and $|h(X_i, \epsilon_i; \mathbf{w}_1) - h(X_i, \epsilon_i; \mathbf{w}_2)| \leq \varrho \|\mathbf{w}_1 - \mathbf{w}_2\|_1$ for any $\mathbf{w}_1, \mathbf{w}_2 \in \mathcal{W}$. Therefore, from Lemma 2,

$$\sup_{\mathbf{w} \in \mathcal{W}} |\text{CV}_1(\mathbf{w})| = O_p(n^{-1/2} |\mathcal{J}|^{1/2}),$$

and we obtain (S3.3).

*Proof of* (S3.4). Take

$$h(X_i, \epsilon_i; \mathbf{w}) = \int_0^{Q_\tau^*(X_i,\mathbf{w}) - Q_\tau(X_i)} \left[ \mathbf{1}\{\varepsilon_i \leq s\} - \mathbf{1}\{\varepsilon_i \leq 0\} \right] ds$$

in Lemma 2. By the fact $|\mathbf{1}\{\varepsilon_i \leq s\} - \mathbf{1}\{\varepsilon_i \leq 0\}| \leq 1$ and Condition 2, we can verify that the conditions of Lemma 2 hold. Therefore, we have $\sup_{\mathbf{w} \in \mathcal{W}} |\text{CV}_2(\mathbf{w})| = O_p(n^{-1/2} |\mathcal{J}|^{1/2})$, and (S3.4).

*Proof of* (S3.5). In view of the fact that $|\psi_\tau(\varepsilon_i)| \leq 1$ and $|\mathbf{1}\{\varepsilon_i \leq s\} - \mathbf{1}\{\varepsilon_i \leq 0\}| \leq 1$, we have

$\sup_{\mathbf{w}\in\mathcal{W}} |CV_{3n}(\mathbf{w})|$

$\leq \frac{2}{n} \sum_{i=1}^{n} \sup_{\mathbf{w}\in\mathcal{W}} \left|\widehat{Q}_\tau(X_i,\mathbf{w}) - Q_\tau^*(X_i,\mathbf{w})\right|$

$\leq 2\max_{1\leq i\leq n} \max_{J_n\in\mathcal{J}} \left|\widehat{Q}_{\tau,J_n}(X_i) - Q_{\tau,J_n}^*(X_i)\right|$

$\leq 2\max_{J_n\in\mathcal{J}} |\widehat{a}_{J_n} - a_{J_n}^*| + 2\max_{1\leq i\leq n}\max_{J_n\in\mathcal{J}} \left|\int_\mathcal{T} \left\{\widehat{b}_{J_n}(t) - b_{J_n}^*(t)\right\}\left\{X_{i,J_n}^*(t) - \mu(t)\right\} dt\right|$

$+ 2\max_{1\leq i\leq n}\max_{J_n\in\mathcal{J}} \left|\int_\mathcal{T} \widehat{b}_{J_n}(t)\left\{\widehat{X}_{i,J_n}(t) - X_{i,J_n}^*(t) - \widehat{\mu}(t) + \mu(t)\right\} dt\right|$

$= O_p(c_n)$. Therefore, (S3.5) is obtained.

*Proof of* (S3.6). By Condition 3, uniformly for $\mathbf{w} \in \mathcal{W}$,

$CV_{4n}(\mathbf{w}) \leq \frac{2}{n} \sum_{k=1}^{K} \sum_{m=1}^{M} \left|\widehat{Q}_\tau(X_{(k-1)M+m},\mathbf{w}) - \widehat{Q}_\tau^{[-k]}(X_{(k-1)M+m},\mathbf{w})\right| = O_p(g_n)$,

which implies (S3.6). This completes the proof of Theorem 1.

## S3.3. Proof of Theorem 2

For any $\mathbf{w}$ (including both random and non-random $\mathbf{w}$'s), we first perform a decomposition of $\widetilde{CV}_K(\mathbf{w})$ as follows $\widetilde{CV}_K(\mathbf{w})$

$= \frac{1}{n}\sum_{k=1}^{K}\sum_{m=1}^{M} \left\{Q_\tau(X_{(k-1)M+m}) - \widehat{Q}_\tau^{[-k]}(X_{(k-1)M+m},\mathbf{w})\right\}\psi_\tau(\varepsilon_{(k-1)M+m})$

$+ \frac{1}{n}\sum_{k=1}^{K}\sum_{m=1}^{M} \int_0^{\widehat{Q}_\tau^{[-k]}(X_{(k-1)M+m},\mathbf{w}) - Q_\tau(X_{(k-1)M+m})} \left[\mathbf{1}\{\varepsilon_{(k-1)M+m} \leq s\} - \mathbf{1}\{\varepsilon_{(k-1)M+m} \leq 0\}\right] ds$

$= \frac{1}{n}\sum_{i=1}^{n} \left\{Q_\tau(X_i) - \widehat{Q}_\tau(X_i,\mathbf{w})\right\}\psi_\tau(\varepsilon_i) + \frac{1}{n}\sum_{i=1}^{n} \int_0^{\widehat{Q}_\tau(X_i,\mathbf{w}) - Q_\tau(X_i)} \left[\mathbf{1}\{\varepsilon_i \leq s\} - \mathbf{1}\{\varepsilon_i \leq 0\}\right] ds$

$+ \frac{1}{n}\sum_{k=1}^{K}\sum_{m=1}^{M} \left\{\widehat{Q}_\tau(X_{(k-1)M+m},\mathbf{w}) - \widehat{Q}_\tau^{[-k]}(X_{(k-1)M+m},\mathbf{w})\right\}\psi_\tau(\varepsilon_{(k-1)M+m})$

$+ \frac{1}{n}\sum_{k=1}^{K}\sum_{m=1}^{M} \int_{\widehat{Q}_\tau(X_{(k-1)M+m},\mathbf{w}) - Q_\tau(X_{(k-1)M+m})}^{\widehat{Q}_\tau^{[-k]}(X_{(k-1)M+m},\mathbf{w}) - Q_\tau(X_{(k-1)M+m})} \left[\mathbf{1}\{\varepsilon_{(k-1)M+m} \leq s\} - \mathbf{1}\{\varepsilon_{(k-1)M+m} \leq 0\}\right] ds$

$= \frac{1}{n}\sum_{i=1}^{n} \left\{Q_\tau(X_i) - \widehat{Q}_\tau(X_i,\mathbf{w})\right\}\psi_\tau(\varepsilon_i) + \frac{1}{n}\sum_{i=1}^{n}\int_0^{\widehat{Q}_\tau(X_i,\mathbf{w}) - Q_\tau(X_i)} \left[\mathbf{1}\{\varepsilon_i \leq s\} - \mathbf{1}\{\varepsilon_i \leq 0\}\right] ds$

$+ O_p(g_n)$,

$= \frac{1}{n}\sum_{i=1}^{n} \left\{Q_\tau(X_i) - \widehat{Q}_\tau(X_i,\mathbf{w})\right\}\psi_\tau(\varepsilon_i) + \frac{1}{n}\sum_{i=1}^{n}\int_0^{\widehat{Q}_\tau(X_i,\mathbf{w}) - Q_\tau(X_i)} \left\{F(s|X_i) - F(0|X_i)\right\} ds$

$+ O_p(g_n) + O_p(c_n) + O_p(n^{-1/2}|\mathcal{J}|^{1/2})$,

where the third equality is due to (S3.2) and the last equality is derived from

$\sup_{\mathbf{w}\in\mathcal{W}} \left|\frac{1}{n}\sum_{i=1}^{n}\int_0^{\widehat{Q}_\tau(X_i,\mathbf{w}) - Q_\tau(X_i)} \left[\mathbf{1}\{\varepsilon_i \leq s\} - \mathbf{1}\{\varepsilon_i \leq 0\} - F(s|X_i) + F(0|X_i)\right] ds\right|$

$\leq \sup_{\mathbf{w}\in\mathcal{W}} \left|\frac{1}{n}\sum_{i=1}^{n}\int_0^{Q_\tau^*(X_i,\mathbf{w}) - Q_\tau(X_i)} \left[\mathbf{1}\{\varepsilon_i \leq s\} - \mathbf{1}\{\varepsilon_i \leq 0\} - F(s|X_i) + F(0|X_i)\right] ds\right|$

$+ \sup_{\mathbf{w}\in\mathcal{W}} \left|\frac{1}{n}\sum_{i=1}^{n}\int_{Q_\tau^*(X_i,\mathbf{w}) - Q_\tau(X_i)}^{\widehat{Q}_\tau(X_i,\mathbf{w}) - Q_\tau(X_i)} \left[\mathbf{1}\{\varepsilon_i \leq s\} - \mathbf{1}\{\varepsilon_i \leq 0\} - F(s|X_i) + F(0|X_i)\right] ds\right|$

$$= O_p(n^{-1/2}|\mathcal{J}|^{1/2}) + O_p(c_n).$$

Here, the last step uses Lemma 2 and a similar argument for the proof of (S3.5). We can see that

$$\widehat{Q}_\tau(X_i, \mathbf{w}) - Q_\tau(X_i)$$
$$= \sum_{J_n \in \mathcal{J}} w_{J_n} \widehat{Q}_{\tau, J_n}(X_i) - Q_\tau(X_i)$$
$$= \sum_{J_n \in \mathcal{J}} w_{J_n} \widehat{a}_{J_n} + \sum_{J_n \in \mathcal{J}} w_{J_n} \int_\mathcal{T} \widehat{b}_{J_n}(t) \{\widehat{X}_{i,J_n}(t) - \widehat{\mu}(t)\} dt - Q_\tau(X_i)$$
$$= \widehat{a}_\mathbf{w} + \sum_{J_n \in \mathcal{J}} w_{J_n} \int_\mathcal{T} \widehat{b}_{J_n}(t) \{\widehat{X}_i(t) - \widehat{\mu}(t)\} dt - Q_\tau(X_i)$$
$$= \widehat{a}_\mathbf{w} + \int_\mathcal{T} \widehat{b}_\mathbf{w}(t) \{\widehat{X}_i(t) - \widehat{\mu}(t)\} dt - a - \int_\mathcal{T} b(t) X_i^c(t) dt$$
$$= \widehat{a}_\mathbf{w} - a + \int_\mathcal{T} \{\widehat{b}_\mathbf{w}(t) - b(t)\} \{\widehat{X}_i(t) - \widehat{\mu}(t)\} dt$$
$$+ \int_\mathcal{T} b(t) \{\mu(t) - \widehat{\mu}(t) + \widehat{X}_i(t) - X_i(t)\} dt. \text{ Thus, we have}$$

$$\frac{1}{n} \sum_{i=1}^n \{\widehat{Q}_\tau(X_i, \mathbf{w}) - Q_\tau(X_i)\} \psi_\tau(\varepsilon_i)$$
$$= (\widehat{a}_\mathbf{w} - a) \frac{1}{n} \sum_{i=1}^n \psi_\tau(\varepsilon_i) + \int_\mathcal{T} \{\widehat{b}_\mathbf{w}(t) - b(t)\} \frac{1}{n} \sum_{i=1}^n \{\widehat{X}_i(t) - \widehat{\mu}(t)\} \psi_\tau(\varepsilon_i) dt + C_1,$$
$$\equiv \Pi_1(\mathbf{w}) + C_1,$$

where $C_1$ is unrelated to $\mathbf{w}$. By Conditions 2 and 5, we have

$$\frac{1}{n} \sum_{i=1}^n \int_0^{\widehat{Q}_\tau(X_i, \mathbf{w}) - Q_\tau(X_i)} \{F(s|X_i) - F(0|X_i)\} ds$$
$$= \frac{1}{n} \sum_{i=1}^n \int_0^{\widehat{Q}_\tau(X_i, \mathbf{w}) - Q_\tau(X_i)} f(s_i^*|X_i) s\, ds$$
$$\leq \frac{\bar{C}}{2n} \sum_{i=1}^n \{\widehat{Q}_\tau(X_i, \mathbf{w}) - Q_\tau(X_i)\}^2$$
$$= \frac{\bar{C}}{2n} \sum_{i=1}^n \left[\widehat{a}_\mathbf{w} - a + \int_\mathcal{T} \{\widehat{b}_\mathbf{w}(t) - b(t)\} \{\widehat{X}_i(t) - \widehat{\mu}(t)\} dt\right]^2$$
$$+ \frac{\bar{C}}{n} \sum_{i=1}^n \left(\left[\widehat{a}_\mathbf{w} - a + \int_\mathcal{T} \{\widehat{b}_\mathbf{w}(t) - b(t)\} \{\widehat{X}_i(t) - \widehat{\mu}(t)\} dt\right]\right.$$
$$\left. \times \int_\mathcal{T} b(t) \{\mu(t) - \widehat{\mu}(t) + \widehat{X}_i(t) - X_i(t)\} dt\right)$$
$$+ \frac{\bar{C}}{2n} \sum_{i=1}^n \left[\int_\mathcal{T} b(t) \{\mu(t) - \widehat{\mu}(t) + \widehat{X}_i(t) - X_i(t)\} dt\right]^2$$
$$\equiv \frac{\bar{C}}{2n} \sum_{i=1}^n \left[\widehat{a}_\mathbf{w} - a + \int_\mathcal{T} \{\widehat{b}_\mathbf{w}(t) - b(t)\} \{\widehat{X}_i(t) - \widehat{\mu}(t)\} dt\right]^2 + \bar{C}\Pi_2(\mathbf{w}) + O_p(d_n^2) + O_p(c_n^2),$$

where $s_i^*$ is between 0 and $\widehat{Q}_\tau(X_i, \mathbf{w}) - Q_\tau(X_i)$ and the last step follows from

$$\max_{1 \leq i \leq n} \left| \int_{\mathcal{T}} b(t) \left\{ \mu(t) - \widehat{\mu}(t) + \widehat{X}_i(t) - X_i(t) \right\} dt \right|$$

$$= \max_{1 \leq i \leq n} \left| \int_{\mathcal{T}} b(t) \left\{ \mu(t) - \widehat{\mu}(t) + \widehat{X}_i(t) - X_{i,J^*}(t) \right\} dt \right|$$

$$\leq \max_{1 \leq i \leq n} \left| \int_{\mathcal{T}} \left\{ b(t) - \widehat{b}_{J^*}(t) \right\} \left\{ \mu(t) - \widehat{\mu}(t) + \widehat{X}_i(t) - X_{i,J^*}(t) \right\} dt \right|$$

$$+ \max_{1 \leq i \leq n} \left| \int_{\mathcal{T}} \widehat{b}_{J^*}(t) \left\{ \mu(t) - \widehat{\mu}(t) + \widehat{X}_i(t) - X_{i,J^*}(t) \right\} dt \right|$$

$$= \max_{1 \leq i \leq n} \left| \int_{\mathcal{T}} \left\{ b(t) - \widehat{b}_{J^*}(t) \right\} \left\{ \mu(t) - \widehat{\mu}(t) + \widehat{X}_i(t) - X_{i,J^*}(t) \right\} dt \right|$$

$$+ \max_{1 \leq i \leq n} \left| \int_{\mathcal{T}} \widehat{b}_{J^*}(t) \left\{ \mu(t) - \widehat{\mu}(t) + \widehat{X}_{i,J^*}(t) - X_{i,J^*}(t) \right\} dt \right|$$

$$= O_p(d_n) + O_p(c_n),$$

which can be easily derived by $J^* \in \mathcal{J}$, (5.6), and Condition 1. In addition,

$$|\Pi_1(\mathbf{w})|$$

$$\leq |\widehat{a}_{\mathbf{w}} - a| \left| \frac{1}{n} \sum_{i=1}^{n} \psi_\tau(\varepsilon_i) \right|$$

$$+ \left[ \int_{\mathcal{T}} \left\{ \widehat{b}_{\mathbf{w}}(t) - b(t) \right\}^2 dt \right]^{1/2} \left( \int_{\mathcal{T}} \left[ \frac{1}{n} \sum_{i=1}^{n} \left\{ \widehat{X}_i(t) - \widehat{\mu}(t) \right\} \psi_\tau(\varepsilon_i) \right]^2 dt \right)^{1/2}$$

$$= \left( |\widehat{a}_{\mathbf{w}} - a| + \left[ \int_{\mathcal{T}} \left\{ \widehat{b}_{\mathbf{w}}(t) - b(t) \right\}^2 dt \right]^{1/2} \right) O_p(n^{-1/2}), \tag{S3.-1}$$

and

$$|\Pi_2(\mathbf{w})|$$

$$\leq O_p(c_n) \left[ |\widehat{a}_{\mathbf{w}} - a| + \frac{1}{n} \sum_{i=1}^{n} \left| \int_{\mathcal{T}} \left\{ \widehat{b}_{\mathbf{w}}(t) - b(t) \right\} \left\{ \widehat{X}_i(t) - \widehat{\mu}(t) \right\} dt \right| \right]$$

$$\leq O_p(c_n) \left( |\widehat{a}_{\mathbf{w}} - a| + \left[ \int_{\mathcal{T}} \left\{ \widehat{b}_{\mathbf{w}}(t) - b(t) \right\}^2 dt \right]^{1/2} \frac{1}{n} \sum_{i=1}^{n} \left[ \int_{\mathcal{T}} \left\{ \widehat{X}_i(t) - \widehat{\mu}(t) \right\}^2 dt \right]^{1/2} \right)$$

$$= \left( |\widehat{a}_{\mathbf{w}} - a| + \left[ \int_{\mathcal{T}} \left\{ \widehat{b}_{\mathbf{w}}(t) - b(t) \right\}^2 dt \right]^{1/2} \right) O_p(c_n). \tag{S3.0}$$

## S3 PROOFS OF THEOREMS 1 AND 2

Let $\mathbf{w}_{\text{true}}$ be a weight vector in which the element corresponding to the true model $J^*$ is 1, and all others are 0's. Using (5.6), (S3.-1), and (S3.0), we have $\Pi_1(\mathbf{w}_{\text{true}}) = O_p(n^{-1/2}d_n)$, $\Pi_2(\mathbf{w}_{\text{true}}) = O_p(d_n c_n)$, and

$$\frac{\bar{C}}{2n} \sum_{i=1}^{n} \left[ \widehat{a}_{\mathbf{w}_{\text{true}}} - a + \int_{\mathcal{T}} \left\{ \widehat{b}_{\mathbf{w}_{\text{true}}}(t) - b(t) \right\} \left\{ \widehat{X}_i(t) - \widehat{\mu}(t) \right\} dt \right]^2$$

$$\leq \bar{C}(\widehat{a}_{\mathbf{w}_{\text{true}}} - a)^2 + \bar{C} \int_{\mathcal{T}} \left\{ \widehat{b}_{\mathbf{w}_{\text{true}}}(t) - b(t) \right\}^2 dt \frac{1}{n} \sum_{i=1}^{n} \int_{\mathcal{T}} \left\{ \widehat{X}_i(t) - \widehat{\mu}(t) \right\}^2 dt$$

$$= O_p(d_n^2). \tag{S3.1}$$

Combining (S3.3)–(S3.1), we know that for $\mathbf{w}_{\text{true}}$,

$$\widetilde{\text{CV}}_K(\mathbf{w}_{\text{true}}) = O_p(d_n^2) + O_p(n^{-1/2}d_n) + O_p(g_n) + O_p(c_n) + O_p(n^{-1/2}|\mathcal{J}|^{1/2}) - C_1. \tag{S3.2}$$

From Condition 5, we also have

$$\frac{1}{n} \sum_{i=1}^{n} \int_0^{\widehat{Q}_\tau(X_i, \mathbf{w}) - Q_\tau(X_i)} \{F(s|X_i) - F(0|X_i)\} \, ds$$

$$\geq \frac{\underline{C}}{2n} \sum_{i=1}^{n} \left[ \widehat{a}_{\mathbf{w}} - a + \int_{\mathcal{T}} \left\{ \widehat{b}_{\mathbf{w}}(t) - b(t) \right\} \left\{ \widehat{X}_i(t) - \widehat{\mu}(t) \right\} dt \right]^2 + \underline{C}\Pi_2(\mathbf{w}) + O_p(d_n^2) + O_p(c_n^2),$$

which along with (S3.3), (S3.3), (S3.2), and $\widetilde{\text{CV}}_K(\widehat{\mathbf{w}}) \leq \widetilde{\text{CV}}_K(\mathbf{w}_{\text{true}})$ implies that

$$\frac{\underline{C}}{2n} \sum_{i=1}^{n} \left[ \widehat{a}_{\widehat{\mathbf{w}}} - a + \int_{\mathcal{T}} \left\{ \widehat{b}_{\widehat{\mathbf{w}}}(t) - b(t) \right\} \left\{ \widehat{X}_i(t) - \widehat{\mu}(t) \right\} dt \right]^2 - \Pi_1(\widehat{\mathbf{w}}) + \underline{C}\Pi_2(\widehat{\mathbf{w}})$$

$$\leq O_p(d_n^2) + O_p(n^{-1/2}d_n) + O_p(g_n) + O_p(c_n) + O_p(n^{-1/2}|\mathcal{J}|^{1/2}). \tag{S3.3}$$

Let $\mathcal{I}$ be the set of $i$ such that the inequality in Condition 7 holds. From Condition 7,

$$\sum_{i=1}^{n}\left[\widehat{a}_{\widehat{\mathbf{w}}} - a + \int_{\mathcal{T}}\left\{\widehat{b}_{\widehat{\mathbf{w}}}(t) - b(t)\right\}\left\{\widehat{X}_i(t) - \widehat{\mu}(t)\right\}dt\right]^2$$

$$\geq \sum_{i\in\mathcal{U}}\lambda_{\min}\left[(\widehat{a}_{\widehat{\mathbf{w}}} - a)^2 + \int_{\mathcal{T}}\left\{\widehat{b}_{\widehat{\mathbf{w}}}(t) - b(t)\right\}^2 dt\right]$$

$$= \lambda_{\min} n^*\left[(\widehat{a}_{\widehat{\mathbf{w}}} - a)^2 + \int_{\mathcal{T}}\left\{\widehat{b}_{\widehat{\mathbf{w}}}(t) - b(t)\right\}^2 dt\right], \tag{S3.4}$$

where $n^*$ is the number of elements in $\mathcal{I}$ and $n^*$ has the same order as $n$ from Condition 7. It follows from (S3.3) and (S3.4) that

$$\frac{Cn^*}{2n}\lambda_{\min}\left[(\widehat{a}_{\widehat{\mathbf{w}}} - a)^2 + \int_{\mathcal{T}}\left\{\widehat{b}_{\widehat{\mathbf{w}}}(t) - b(t)\right\}^2 dt\right]$$

$$\leq O_p(d_n^2) + O_p(n^{-1/2}d_n) + O_p(g_n) + O_p(c_n) + O_p(n^{-1/2}|\mathcal{J}|^{1/2}) + \Pi_1(\widehat{\mathbf{w}}) - \underline{C}\Pi_2(\widehat{\mathbf{w}})$$

$$\leq \left(|\widehat{a}_{\widehat{\mathbf{w}}} - a| + \left[\int_{\mathcal{T}}\left\{\widehat{b}_{\widehat{\mathbf{w}}}(t) - b(t)\right\}^2 dt\right]^{1/2}\right)\{O_p(n^{-1/2}) + O_p(c_n)\}$$

$$+ O_p(d_n^2) + O_p(n^{-1/2}d_n) + O_p(g_n) + O_p(c_n) + O_p(n^{-1/2}|\mathcal{J}|^{1/2}).$$

Thus, there exists $\widetilde{c}_n = O_p(n^{-1/2}) + O_p(c_n)$ such that

$$\frac{Cn^*}{2n}\lambda_{\min}\left[(\widehat{a}_{\widehat{\mathbf{w}}} - a)^2 + \int_{\mathcal{T}}\left\{\widehat{b}_{\widehat{\mathbf{w}}}(t) - b(t)\right\}^2 dt\right] + \widetilde{c}_n\left(|\widehat{a}_{\widehat{\mathbf{w}}} - a| + \left[\int_{\mathcal{T}}\left\{\widehat{b}_{\widehat{\mathbf{w}}}(t) - b(t)\right\}^2 dt\right]^{1/2}\right)$$

$$\leq O_p(d_n^2) + O_p(n^{-1/2}d_n) + O_p(g_n) + O_p(c_n) + O_p(n^{-1/2}|\mathcal{J}|^{1/2}).$$

This leads to

$$(|\widehat{a}_{\widehat{\mathbf{w}}} - a| + \underline{C}^{-1}\lambda_{\min}^{-1}\widetilde{c}_n n/n^*)^2 + \left(\left[\int_{\mathcal{T}}\left\{\widehat{b}_{\widehat{\mathbf{w}}}(t) - b(t)\right\}^2 dt\right]^{1/2} + \underline{C}^{-1}\lambda_{\min}^{-1}\widetilde{c}_n n/n^*\right)^2$$

$$\leq O_p(d_n^2) + O_p(n^{-1/2}d_n) + O_p(g_n) + O_p(c_n) + O_p(n^{-1/2}|\mathcal{J}|^{1/2}) - 2\underline{C}^{-2}\lambda_{\min}^{-2}\widetilde{c}_n^2 n^2/n^{*2},$$

which, along with $\widetilde{c}_n = O_p(n^{-1/2}) + O_p(c_n)$, yields that

$$|\widehat{a}_{\widehat{\mathbf{w}}} - a| + \left[\int_{\mathcal{T}} \left\{\widehat{b}_{\widehat{\mathbf{w}}}(t) - b(t)\right\}^2 dt\right]^{1/2} + 2\underline{C}^{-1}\lambda_{\min}^{-1}\widetilde{c}_n n/n^*$$

$$= O_p(d_n) + O_p(n^{-1/4}d_n^{1/2}) + O_p(g_n^{1/2}) + O_p(c_n^{1/2}) + O_p(n^{-1/4}|\mathcal{J}|^{1/4}).$$

It further leads to

$$|\widehat{a}_{\widehat{\mathbf{w}}} - a| + \left[\int_{\mathcal{T}} \left\{\widehat{b}_{\widehat{\mathbf{w}}}(t) - b(t)\right\}^2 dt\right]^{1/2}$$

$$= O_p(d_n) + O_p(n^{-1/4}d_n^{1/2}) + O_p(g_n^{1/2}) + O_p(c_n^{1/2}) + O_p(n^{-1/4}|\mathcal{J}|^{1/4}),$$

and thus

$$(\widehat{a}_{\widehat{\mathbf{w}}} - a)^2 + \int_{\mathcal{T}} \left\{\widehat{b}_{\widehat{\mathbf{w}}}(t) - b(t)\right\}^2 dt = O_p(d_n^2 + g_n + c_n + n^{-1/2}|\mathcal{J}|^{1/2}).$$

This completes the proof of Theorem 2.

## S4.   Simulation Studies Continued From Section 6

### S4.1.   Simulation Setting and Additional Figures for Simulation Design I

We assume that the functional covariate $X_i$ is observed at time points $T_{il}$ with measurement errors for $l = 1, \ldots, N_i$. We consider the following simulation designs.

(a) Consider sparse designs for $X_i$. Set $n_T = n + n_0$ as the sample size of a dataset, among which the first $n$ dataset are used as training data and the other $n_0 = 100$ as test data. Set $N_i$ are sampled from the discrete uniform distribution on $\{10, 11, 12\}$. For each $i$, $T_{il}$, $l = 1, \ldots, N_i$ are sampled from a uniform distribution on $\mathcal{T} = [0, 1]$.

(b) Set the eigenfunctions $\phi_j(t) = \sqrt{2}\cos(j\pi t)$ and the eigenvalues $\kappa_j = j^{-1.2}$, $j = 1, 2, \ldots$. Set $\mu(t) = 0$ and generate $X_i(t) = \mu(t) + \sum_{j=1}^{J_{\max}} \kappa_j^{1/2} Z_{ij}\phi_j(t)$, where $J_{\max} = 20$ and $Z_{ij}$ are independently sampled from a uniform distribution on $[-3^{1/2}, 3^{1/2}]$ for each $j$.

(c) Obtain the observations $U_{il}$ by adding measurement errors $\epsilon_{il}$ to $X_i(T_{il})$, i.e., $U_{il} = X_i(T_{il}) + \epsilon_{il}$, where $\epsilon_{il}$ are sampled from $N(0, 0.8)$.

(d) Set $b^{[1]}(t) = \sum_{j=1}^{J_{\max}} j^{-1}\phi_j(t)$, $a^{[2]} = 2\sum_{j=1}^{J_{\max}} j^{-1.5}\kappa_j^{1/2}$, and $b^{[2]}(t) = \sum_{j=1}^{J_{\max}} j^{-1.5}\phi_j(t)$. Generate the response observation $Y_i$ by the following heteroscedastic functional linear model

$$Y_i = \theta \int_0^1 b^{[1]}(t) X_i(t)\, dt + \sigma(X_i) e_i,$$

where $\sigma(X_i) = a^{[2]} + \int_0^1 b^{[2]}(t) X_i(t)\, dt$ and $e_i$ are sampled from $N(0,1)$. It is easy to see that $\sigma(X_i) > 0$ for $X_i$ in its domain.

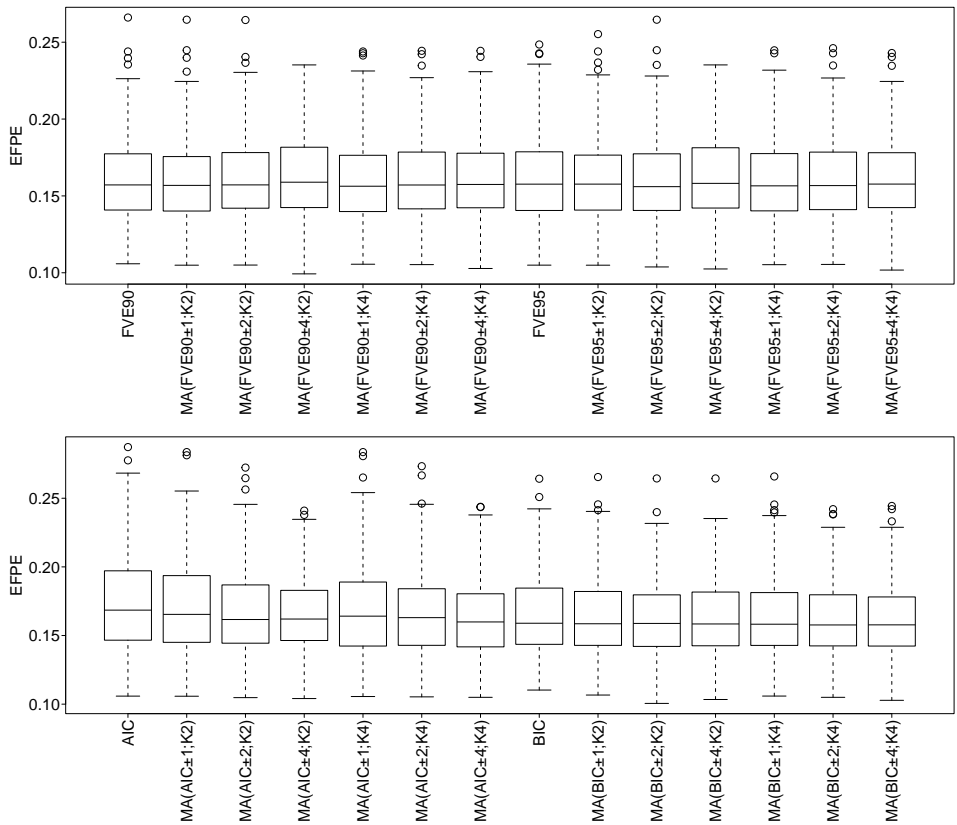

**Figure S.9** Boxplots of EFPE in simulation design I with $\tau = 0.5$. FVE90, FVE with $\gamma = 0.90$; FVE95, FVE with $\gamma = 0.95$; MA(FVE90±$\alpha$, K$\beta$), model averaging method with $\widehat{J}_n$ determined by FVE90, with $d = \alpha$, and with weights selected by $\mathrm{CV}_{K=\beta}$; MA(FVE95±$\alpha$, K$\beta$), MA(AIC±$\alpha$, K$\beta$), and MA(BIC±$\alpha$, K$\beta$) have similar definitions.

## S4.2. Simulation Design II

We consider the same simulation design as that in Subsection S4.1 except for $J_{\max}$, $b^{[1]}(t)$, $a^{[2]}$, and $b^{[2]}(t)$. Specifically, we set $J_{\max} = 8$, $b^{[1]}(t) = \sum_{j=1}^{3} j^{-1}\phi_j(t)$, $a^{[2]} = 2\sum_{j=1}^{3} j^{-1.5}\kappa_j^{1/2}$, and $b^{[2]}(t) = \sum_{j=1}^{3} j^{-1.5} \phi_j(t)$. The candidate set $\mathcal{J}$ is fixed to be $\{0, 1, \ldots, 6\}$. Accordingly, the true models are $J^* = 3, 4, 5,$ or $6$, which are included in the candidate models with the candidate set $\mathcal{J}$.

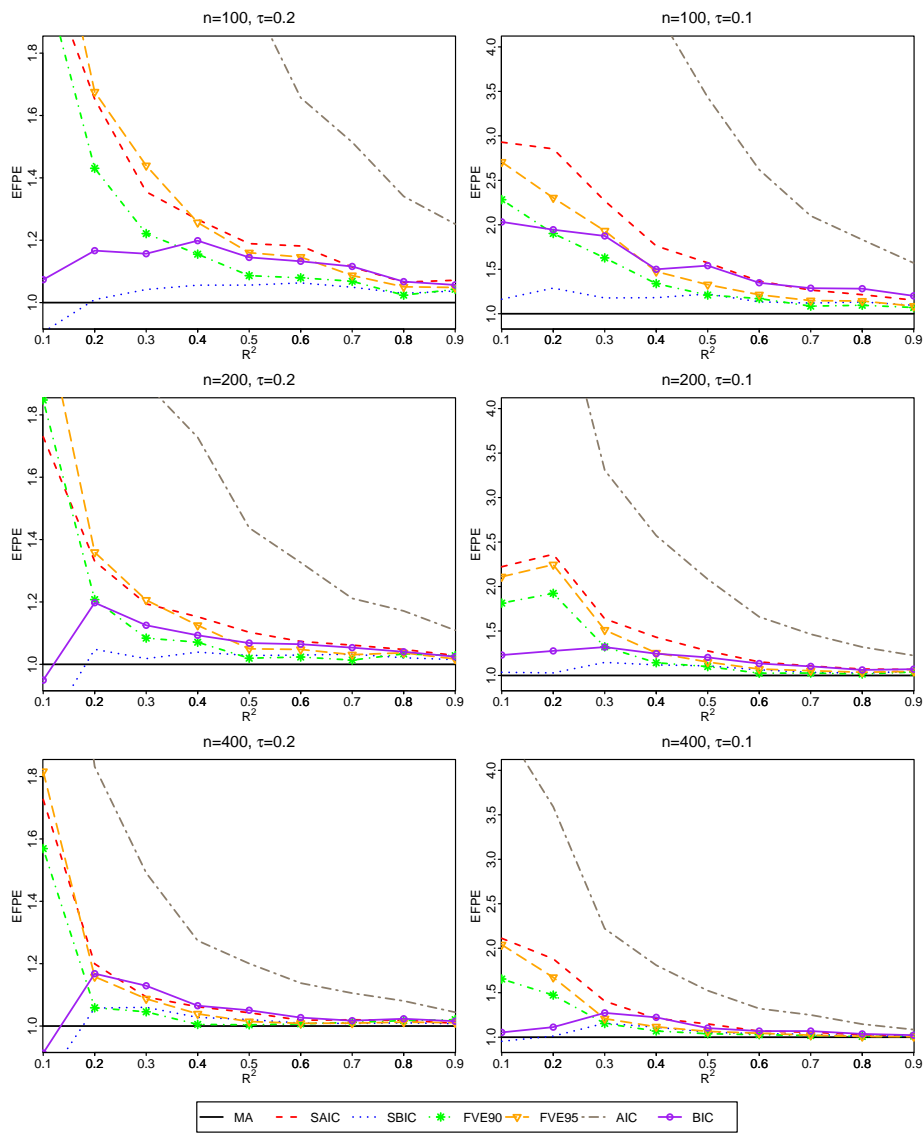

**Figure S.10** Normalized EFPE in simulation design I for $\tau = 0.2$ and $0.1$. MA denotes MA(FVE90$\pm 4$, $K4$).

We first evaluate the finite sample performance of the proposed model averaging prediction in the scenario with true candidate models. As in Subsection S4.1, we compare it with the other six model averaging and selection methods by computing the excess final prediction errors. We set $K = 4$ and consider $n = 100, 200, 400$ and $\tau = 0.5, 0.05$. In each simulation setting, we normalize the EFPEs of the other six methods by a division of the EFPE of our method. The results are presented in Figure

S.11. The performance of different methods becomes more similar when $R^2$ is large. When $\tau = 0.5$, the performance of our method is the best in some cases (e.g., $R^2 = 0.1, 0.4, 0.5$) for $n = 100$, while our method has no advantage for $n = 200, 400$. When $\tau = 0.05$, our method significantly outperforms the other six methods when $R^2 \leq 0.6$. To further illustrate the advantage of our method in terms of prediction for extreme quantiles, we conduct additional simulation studies with $\tau = 0.2$ and $0.1$. The results are displayed in Figure S.12. From Figures S.11 and S.12, the advantage of our method over the other six methods increases as $\tau$ decreases. This illustrates that even when there is at least one true model in the set of candidate models, our method leads to smaller EFPEs compared with the other six methods for extreme quantiles.

Next, we verify the consistency of $\widehat{b}_{\widehat{\mathbf{w}}}(t)$ in Theorem 2. We do this by computing the mean integrated squared error (MISE) based on 200 replications. Specifically, the MISE of our estimator is computed as

$$\text{MISE} = \frac{1}{200} \sum_{r=1}^{200} \int_0^1 \left\{ \widehat{b}_{\widehat{\mathbf{w}}}(t)^{(r)} - b(t) \right\}^2 dt,$$

where $\widehat{b}_{\widehat{\mathbf{w}}}(t)^{(r)}$ denotes the model averaging estimator of $b(t)$ in the $r$th replication. We set $K = 4$ and consider $n = 50, 100, 300, 500, 700, 900, 1100$, $R^2 = 0.1, 0.5, 0.9$, and $\tau = 0.5, 0.05$. Figure S.13 displays MISE against $n$ for each combination of $R^2$ and $\tau$. These plots show that the MISE of $\widehat{b}_{\widehat{\mathbf{w}}}(t)$ decreases to zero as $n$ increases, which reflects the consistency of $\widehat{b}_{\widehat{\mathbf{w}}}(t)$. In addition, we observe that the MISE increases as $R^2$ increases, which is not counterintuitive since larger values of $R^2$ correspond to larger $\theta$ and $b(t)$.

Finally, we compare the performance of $\widehat{b}_{\widehat{\mathbf{w}}}(t)$ with the other six methods. We set $K = 4$ and consider $n = 100, 200, 400$ and $\tau = 0.5, 0.05, 0.01$. In each simulation setting, we normalize the MISEs of the other six methods by dividing them by the MISE of our method. The results are presented in Figure S.14. When $\tau = 0.5$, the performance of FVE(0.90) is the best, AIC and SAIC are the worst, and our method has no advantage over BIC, SBIC, and FVE(0.95). When $\tau = 0.05$, our method and FVE(0.90) outperform the other five methods. Further, when $\tau$ decreases to $0.01$, our method has by far the best performance. Overall, the results illustrate the advantage of our method for extreme quantiles.

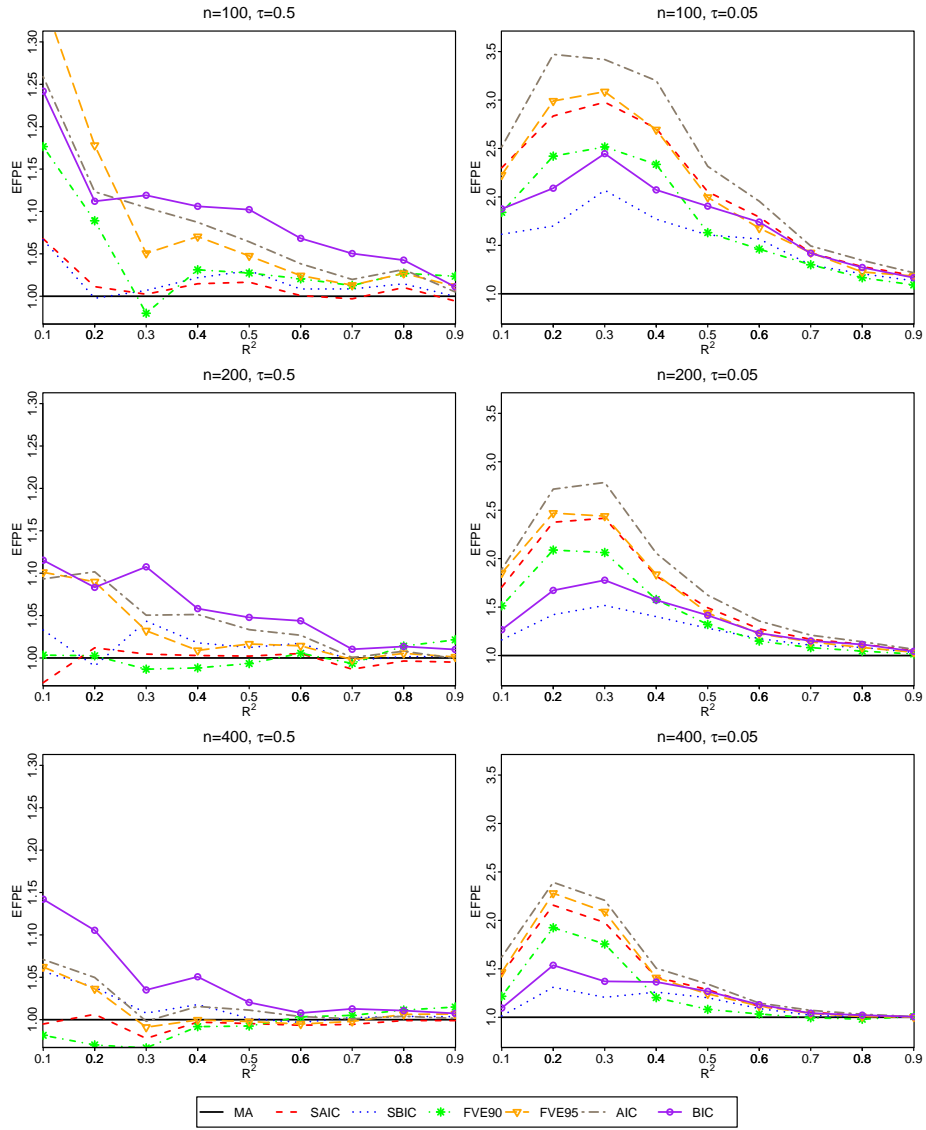

**Figure S.11**    Normalized EFPE in simulation design II for $\tau = 0.5$ and $0.05$. MA denotes MA(FVE90$\pm 4$, $K4$).

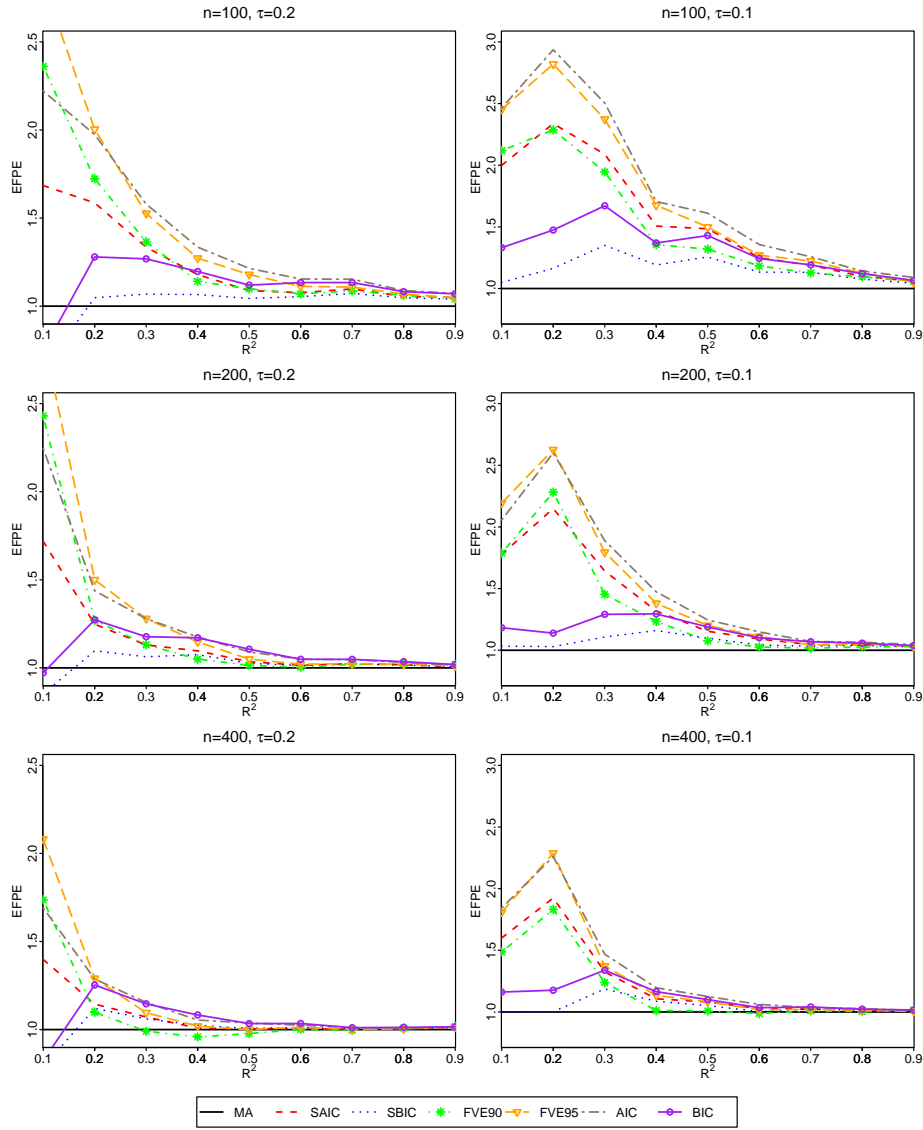

**Figure S.12** Normalized EFPE in simulation design II for $\tau = 0.2$ and $0.1$. MA denotes MA(FVE90$\pm$4, $K$4).

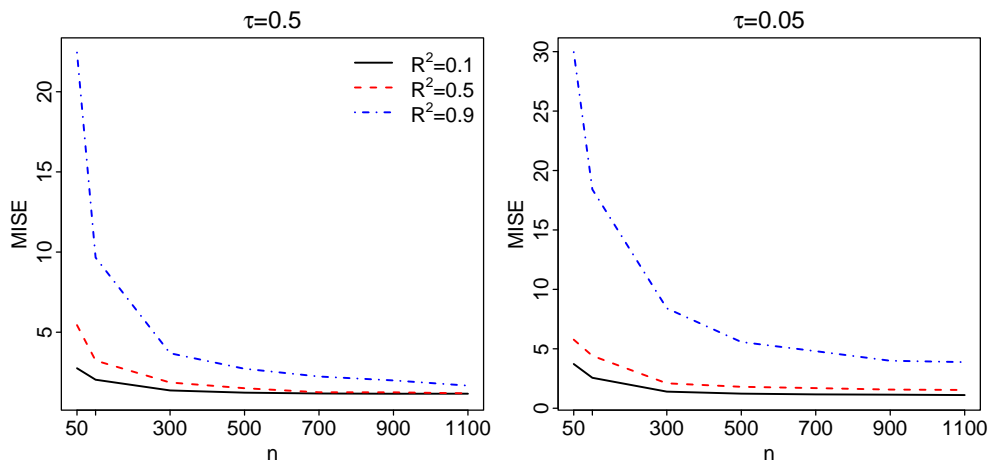

**Figure S.13**   MISE of $\widehat{b}_{\widehat{\mathbf{w}}}(t)$ with $\tau = 0.5$ **(left)** and $\tau = 0.05$ **(right)** under simulation design II.

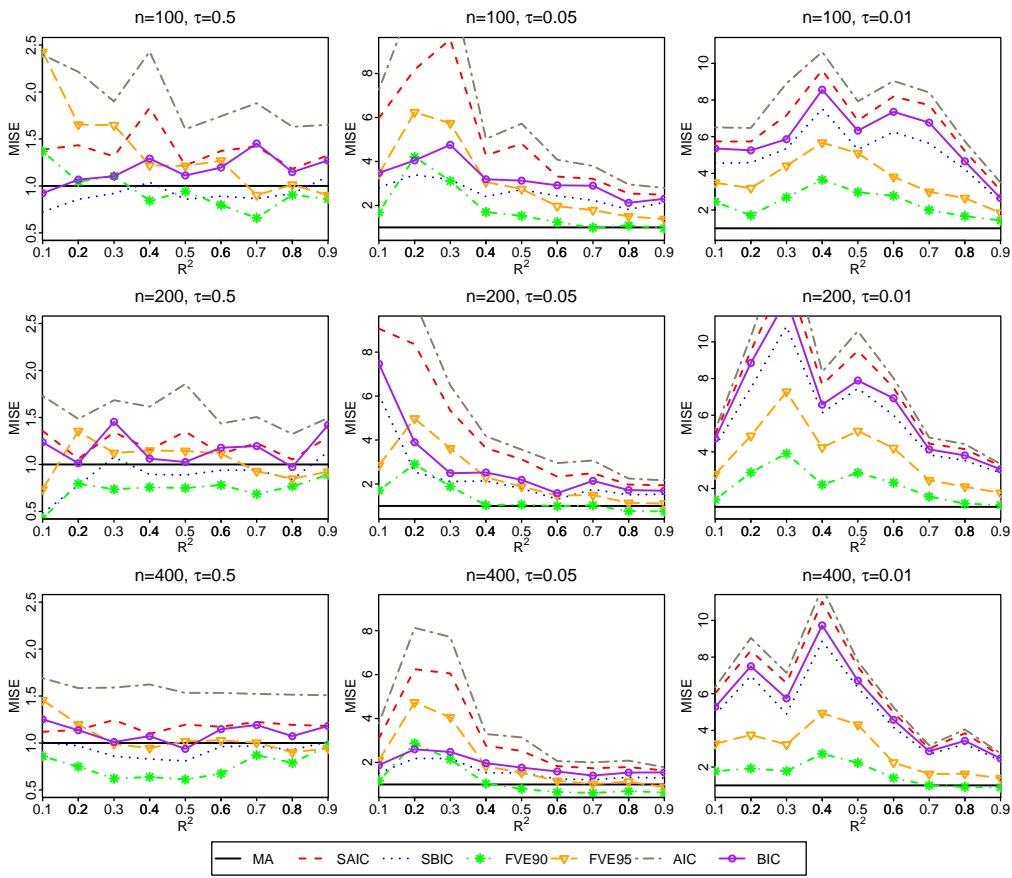

**Figure S.14**   Normalized MISE in simulation design II for $\tau = 0.5$ **(left)**, $0.05$ **(middle)**, and $0.01$ **(right)**. MA denotes **MA(FVE90±4, $K$4)**.